\documentclass{article}

\usepackage{PRIMEarxiv}

\usepackage[utf8]{inputenc} % allow utf-8 input
\usepackage[T1]{fontenc}    % use 8-bit T1 fonts
\usepackage{hyperref}       % hyperlinks
\usepackage{natbib}
\usepackage{doi}
\usepackage{url}            % simple URL typesetting
\usepackage{booktabs}       % professional-quality tables
\usepackage{amsfonts}       % blackboard math symbols
\usepackage{amsmath}
\usepackage{amssymb}
\usepackage{float}
\usepackage{amsthm}
\usepackage{nicefrac}       % compact symbols for 1/2, etc.
\usepackage{microtype}      % microtypography
\usepackage{lipsum}
\usepackage{fancyhdr}       % header
\usepackage{graphicx}       % graphics

\usepackage{enumerate}

\RequirePackage{tikz}
\RequirePackage{tikz-qtree}

\theoremstyle{thmstyleone}%
\newtheorem{theorem}{Theorem}%  meant for continuous numbers
\newtheorem{proposition}[theorem]{Proposition}% 
\newtheorem{lemma}[theorem]{Lemma}
\newtheorem{corollary}[theorem]{Corollary}

\theoremstyle{thmstyletwo}%
\newtheorem{example}{Example}%
\newtheorem{remark}{Remark}%

\theoremstyle{thmstylethree}%
\newtheorem{definition}{Definition}%

\graphicspath{{media/}}     % organize your images and other figures under media/ folder
\usepackage{commands}

\newcommand{\REV}[2]{#2}

\title{Finite Hilbert systems for Weak Kleene logics
}

\author{\href{https://orcid.org/0000-0003-3240-386X}{\includegraphics[scale=0.06]{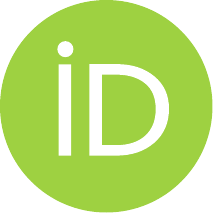}}\hspace{1mm}%
  Vitor Greati \\
  Bernoulli Institute \\
  University of Groningen \\
  Groningen, The Netherlands\\
  \texttt{v.rodrigues.greati@rug.nl} \\
   \And
  \href{https://orcid.org/0000-0002-6941-7555}{\includegraphics[scale=0.06]{orcid.pdf}\hspace{1mm}}Sérgio Marcelino \\
  SQIG – Instituto de Telecomunicações\\
  Departamento de Matemática -- Instituto Superior Técnico \\
  Lisboa, Portugal\\
\texttt{smarcel@math.tecnico.ulisboa.pt} \\
  \And
  \href{https://orcid.org/0000-0003-1364-5003}{\includegraphics[scale=0.06]{orcid.pdf}}\hspace{1mm}%
  Umberto Rivieccio \\
  Departamento de Lógica, Historia y Filosofía de la Ciencia \\
  Universidad Nacional de Educación a Distancia \\
  Madrid, Spain\\
\texttt{umberto@fsof.uned.es} \\
}

\begin{document}
\maketitle

\begin{abstract}
% In this paper we present the first 
% finite Hilbert-style single-conclusion axiomatization of Paraconsistent Weak Kleene and
% Bochvar-Kleene logics. 
% These axiomatizations are obtained by tweaking known multiple-conclusion axiomatizations,
% whose more expressive formalism allows to finitely axiomatize every logic characterized by finite matrix.
%
  % Finally, by carefully %m
  % replacing the multiple-conclusion rules by single-conclusion ones,
  % we obtain the first finite Hilbert-style single-conclusion axiomatizations for these logics.
%
% We present the first finite Hilbert-style single-conclusion axiomatization of Paraconsistent Weak Kleene and Bochvar-Kleene logics. These axiomatizations are obtained by modifying known multiple-conclusion axiomatizations, which, due to their more expressive formalism, are known to enable finite (and analytical) axiomatizations for every logic characterized by a finite matrix."
%
%   Multiple-conclusion Hilbert-style systems are known to enable finite % (and analytical) 
%   axiomatizations for every logic characterized by a finite matrix.
%   We first obtain such axiomatizations for Paraconsistent Weak Kleene and Bochvar-Kleene logics which are both $3$-valued.
% Finally, we modify these axiomatizations, by replacing the multiple-conclusion rules by carefully picked single-conclusion ones, we obtain the first finite Hilbert-style single-conclusion axiomatizations for these logics.
Multiple-conclusion Hilbert-style systems allow us to
finitely axiomatize
%to produce a finite axiomatization for 
every logic defined by a finite matrix. %In this paper we % first 
%obtain 
Having obtained such axiomatizations  for Paraconsistent Weak Kleene and Bochvar-Kleene logics, %, both of which are 3-valued. Finally, 
%and %We 
we modify them %axiomatizations 
by replacing the multiple-conclusion rules with carefully selected single-conclusion ones. 
%resulting in 
In this way we manage to introduce
the first \emph{finite} Hilbert-style single-conclusion axiomatizations for these logics.
\end{abstract}

% keywords can be removed

\section{Introduction}\label{sec:introduction}
%\begin{itemize}

%\item Motivations of BK and PWK\\

In his classic book~\cite{Kl50},
S.C.~Kleene employs two different sets of three-valued truth tables
to introduce the logical systems
 known, in today's parlance, as \emph{Strong Kleene} and \emph{Weak Kleene} logics.
 The latter,
 %which was 
 independently considered in 1937 by Bochvar~\cite{bochvar1938original,bochvar1981},
 is also called  \emph{Bochvar-Kleene  logic} (henceforth $\BKName$).

From a formal point of view, the main difference between the strong  and the weak Kleene tables
is that in the latter the third truth value ($\uv$)
exhibits an \emph{infectious}
behaviour: any interaction 
between  $\uv$ and either of  the  classical values ($\tv$ and $\fv$)
  delivers
 $\uv$ itself. This feature, %of the Bochvar-Kleene tables, 
 as we shall see, makes the resulting 
 logics somewhat less tractable than most well-known three-valued logics, both from an algebraic and a proof-theoretic point of view. 

 From the Bochvar-Kleene tables two logics naturally arise. One ($\BKName$)
is obtained by choosing the single truth value $\tv$ as designated; the other,
which we call
\emph{Paraconsistent Weak Kleene} (\PWKName),
results from designating both $\tv$ and $\uv$. 
Concerning both  these systems, a positive and a negative result are particularly 
worth mentioning in the present context. 

The good news is that
both logics are closely related, from a formal point of view, to the classical: more precisely,
$\PWKName$ and $\BKName$ are, respectively, the \emph{left} and \emph{right variable inclusion companions} of classical logic (more on this below). The bad news, on the other hand, 
is that 
for neither of these logics a \emph{finite} Hilbert-style axiomatization currently exists (this observation was made in~\cite{3valbook}, to which we also refer the reader for further background 
and examples of axiomatizations of three-valued logics). 
In other words, we do not know
whether these logics admit a \emph{finite basis}~\cite{W88}.
This is precisely the gap  we wish to bridge in the present paper, thus solving a fundamental
open problem concerning these logics.

According to Bochvar's original paper~\cite{bochvar1938original,bochvar1981},
the intended applications
of %$\BKName$
the
Bochvar-Kleene  logic %($\BKName$) was introduced in 1937 by Bochvar~\cite{bochvar1981}, 
%who had in mind 
%potential 
%applications 
are in the formalization of paradoxes, future contingent statements
and presuppositions (see e.g.~\cite{ferguson2014}
%A computational interpretation of conceptivism TM Ferguson 
for a more recent computational interpretation of \BKName{}).
The third value is therefore meant to
represent nonsensical statements,
or corrupted data in the  interpretation given by Kleene.
This explains the infectious behaviour of the third value,
because any complex formula having a nonsensical or paradoxical subformula 
should be regarded as nonsensical/paradoxical too.
%and has an \emph{infectious} behaviour: any interaction  of $\uv$ with the other truth values delivers $\uv$ itself.

%The 
Paraconsistent Weak Kleene (\PWKName) seems to have been  considered already by S.~Halld\'en
in his 1949 monograph~\cite{hallden1949}, %Logic of nonsense
and two decades later by A.~Prior~\cite{prior1967}, %Prior Arthur (1967) Past, Present and Future, Oxford, Oxford University Press.
but has only recently been studied in more depth (see e.g.~\cite{BonGilPaoPer17} 
%Ciuni, Carrara...al the papers by Pra Baldi etc. 
and~\cite{pailos2018}, the latter of which
explores applications to the theory of truth).
%Theories of truth based on four-valued infectious logicsBruno Da R ́e1,2, Federico Pailos1,2, and Damian Szmuc1,2

% CHECK MORE MOTIVATION in 3-valued paper!
 %\item 
    %\item  

    \smallskip

    \REV{1.1}{The proof theory of \BKName{} and \PWKName{} 
    has been intensively developed in the
    last years employing different
    formalisms and approaches, like
    sequent calculi~\cite{Paoli2020,Bonzio2022}, natural deduction~\cite{Petrukhin2017,belikov2021}
    and tableaux~\cite{Paoli2020,Bonzio2022}.
    %A number of Hilbert-like proof systems for \BKName{} and \PWKName{} 
    %exist in the literature. 
    In addition, a number of Hilbert-like systems for these logics exist in the literature~\cite{baaz1996,bonzio2021plonka,Bonzio2021-BONCLA}.
    However,
    as explained in
      Section~\ref{sec:Hilbertcal},
    none of them \REV{2.1}{are} finite Hilbert-style systems in the usual sense (we shall call these \emph{\SetFmla{} H-systems}).}
     %Many axiomatizations are known \cite{??somewithsequents,andtableuax?,naturaldeuction}, however when we restrict to Hilbert calculi (here dubbed $\SetFmla{}$ H-systems) no known finite axiomatization for these logics was previously known.
    %
%     Existing Hilbert-style systems:

For 
\BKName{},  a {finite} but non-standard axiomatization may be obtained
by taking
  any complete \SetFmla{} H-system
  %\SetFmla{} H-calculus 
  for classical logic (with \emph{modus ponens} as its only rule) and, while keeping all the axioms, replacing \emph{modus ponens} by  a restricted version that satisfies the \emph{containment condition}~\cite[Prop.~4]{BonGilPaoPer17},~\cite{bonzio2021plonka}. 
  The finite
  Hilbert-style system for \BKName{}
  we introduce here %, in contrast, 
  will instead be standard, i.e.~consisting of a finite number of axioms and unrestricted rule schemas.
  %In the present paper we will aim at introducing a standard
  %Hilbert-style system for \BKName{}, consisting of a finite number of axiom and unrestricted rule schemata.

For both \BKName{} and  \PWKName{}, \emph{infinite} Hilbert-style systems
may be found in~\cite{bonzio2021plonka,Bonzio2021-BONCLA}; we note that the completeness proofs
found 
%introduced in the latter 
in these papers are essentially algebraic, and rely on
the above-mentioned observation that \BKName{} and 
\PWKName{} are, respectively,
the right and the left variable inclusion companion
of classical logic~\cite[Thm.~4, p.~258]{urquhart2001},~\cite{bonzio2021plonka,caleiro2020infec}. 
%and
%the logic \BKName{} dealt with in the sequel is the right variable inclusion companion of classical logic \cite{ISMVL}.
%Logics of left variable inclusion and Płonka sums of matrices S. Bonzio, T. Moraschini & M. Pra Baldi
%\BKName{}
%is the \emph{right variable inclusion companion} of classical logic~\cite[Thm.~4, p.~258]{urquhart2001},

%Using the facts that \PWKName{} is the \emph{left variable inclusion companion}
%of classical logic~\cite{bonzio2021plonka,ISMVL} and
%the logic \BKName{} dealt with in the sequel is the right variable inclusion companion of classical logic \cite{ISMVL}.
%Logics of left variable inclusion and Płonka sums of matrices S. Bonzio, T. Moraschini & M. Pra Baldi
%\BKName{}
%is the \emph{right variable inclusion companion} of classical logic~\cite[Thm.~4, p.~258]{urquhart2001},
%it was shown using algebraic techniques infinite 
%\SetFmla{} axiomatizations for \BKName{} and 
%\PWKName{} where obtained in~\cite{bonzio2021,Bonzio2021-BONCLA}. 

In the present paper  we  follow a two-step strategy.
Relying on the general observation
 that every finite logical matrix can be finitely axiomatized by means of  a Hibert-style
    multiple-conclusion system (here called a \emph{$\SetSet{}$ H-system}),
    we  first introduce finite $\SetSet{}$ H-systems for \BKName{} and  \PWKName{},
     then show how from these
    $\SetFmla{}$ axiomatizations may be obtained preserving finiteness.

The paper is organized as follows. In Section~\ref{sec:preliminaries} we formally introduce
the language and semantics of \BKName{} and  \PWKName{}. Section~\ref{sec:Hilbertcal}
contains as much theory of $\SetSet{}$ and $\SetFmla{}$ H-systems as we shall need in order to introduce our
axiomatic systems for \PWKName{} and   \BKName{}. The former  is then presented and shown to be complete 
in Section~\ref{sec:pwk} (\PWKName{}),
 the latter  in Section~\ref{sec:bk} (\BKName{}).
The final Section~\ref{sec:final_remarks} contains  concluding remarks and suggestions for 
future research.

\section{Language and semantics of $\BKName{}$ and  $\PWKName{}$}\label{sec:preliminaries}
%Let $\SigInfec \SymbDef \Set{\land,\lor,\neg}$ be a propositional signature, such that $\land$
%and 
Let $\land$,
$\lor$ and $\to$ be binary connectives
and $\neg$ be a unary connective.
Call a collection $\SigA$ of these connectives 
a \emph{propositional signature}.
We may write
$\SigA_{\ConA_1\ldots\ConA_n}$ for the signature
$\Set{\ConA_1,\ldots,\ConA_n} \subseteq \Set{\land,\lor,\to,\neg}$.

A \emph{$\SigA$-algebra} is
a structure $\AlgA \SymbDef \Struct{\ValSetA, \AlgInterp{\cdot}{\AlgA}}$
such that $\ValSetA$ is a 
nonempty set called the \emph{carrier}
of $\AlgA$ and, for each
$k$-ary connective $\ConA \in \SigA$,
the $k$-ary mapping
$\AlgInterp{\ConA}{\AlgA} : \ValSetA^k \to \ValSetA$ is the \emph{interpretation} (or \emph{truth table}) \emph{of $\ConA$ in $\AlgA$}.

Given a denumerable
set $\PropSet$ of \emph{propositional variables},
we denote by $\LangAlg{\SigA}{\PropSet}$
%the absolutely free algebra 
the term algebra
over $\SigA$ 
%freely
generated by $\PropSet$ or, more briefly, the \emph{$\SigA$-language (generated by $\PropSet$)}, whose universe
is denoted by $\LangSet{\SigA}{\PropSet}$.
The elements of the latter are called \emph{$\SigA$-formulas}.
Propositional variables will be denoted
by lowercase letters $\PropA,\PropB,\PropC,\PropD$,
and $\SigA$-formulas will be denoted by Greek letters
$\FmA,\FmB,\FmC,\FmD$,
possibly subscripted with positive integers.

%We  define the binary connective $\imp$
%by abbreviation: for all $\FmA,\FmB \in \LangSetInfec$,
%let $\FmA\imp\FmB \SymbDef \neg\FmA\lor\FmB$.
The endomorphisms on $\LangAlg{\SigA}{\PropSet}$
are called \emph{$\SigA$-substitutions}.
By $\Subf(\FmSetA)$ we denote the set of all subformulas of the formulas in $\FmSetA \subseteq \LangSet{\SigA}{\PropSet}$.
Moreover, we will usually 
write $\FmSetA,\FmSetB$
to denote $\FmSetA \cup \FmSetB$
and we will omit curly braces when writing
sets of formulas.
\REV{1.3,2.3}{Also, we write
$\FmSetCompl{\FmSetA}$
for $\LangSet{\SigA}{\PropSet}{\setminus}\FmSetA$}.

We take $\SigInfec$ to be the signature
of classical logic as well as that of $\PWKName$ and $\BKName$ in the present work.
We are going to define these logics in a moment via matrix semantics.

Let $\BoolAlg \SymbDef \Struct{\Set{\fv,\tv}, \AlgInterp{\cdot}{\BoolAlg}}$
be the standard \REV{1.2}{two-element} Boolean $\SigInfec$-algebra. 
For
$\ValSetInfec \SymbDef \Set{\fv,\uv,\tv}$,
define the $\SigInfec$-algebra $\AlgInfec \SymbDef \Struct{\ValSetInfec, \AlgInterp{\cdot}{\AlgInfec}}$
such that the connectives in $\SigInfec$
are interpreted according to the following
truth tables:

\begin{center}
    \vspace{1em}
    \begin{tabular}{c|ccc}
         $\AlgInterp{\land}{\AlgInfec}$& \fv & \uv & \tv \\
         \midrule
         \fv&\fv&\uv&\fv\\
         \uv&\uv&\uv&\uv\\
         \tv&\fv&\uv&\tv\\
    \end{tabular}\quad
    \begin{tabular}{c|ccc}
         $\AlgInterp{\lor}{\AlgInfec}$& \fv & \uv & \tv \\
         \midrule
         \fv&\fv&\uv&\tv\\
         \uv&\uv&\uv&\uv\\
         \tv&\tv&\uv&\tv\\
    \end{tabular}\quad
    \begin{tabular}{c|ccc}
         & $\AlgInterp{\neg}{\AlgInfec}$\\
         \midrule
         \fv&\tv\\
         \uv&\uv\\
         \tv&\fv\\
    \end{tabular}
    \vspace{1em}
\end{center}

\noindent
As we will see in a moment, such $\SigInfec$-algebra provides
the interpretation structure
for the logical matrices that determine the logics $\PWKName{}$ and $\BKName{}$.
Note that
we have,
for all $\ConA \in \SigInfec$ of arity $k$,
$\AlgInterp{\ConA}{\AlgInfec}(\vec{\ValA})
= \AlgInterp{\ConA}{\BoolAlg}(\vec{\ValA})$
if $\vec{\ValA} \in \Set{\fv,\tv}^k$
and $\AlgInterp{\ConA}{\AlgInfec}(\vec{\ValA})
= \uv$ otherwise.
In other words, the above truth tables result from
extending the classical two-valued  tables
with an \emph{infectious truth value}~\cite{caleiro2020infec}.

We now extend the above observation
to \REV{2.2}{the} derived operations of $\AlgInfec$.
%Before doing so, let us introduce some
%notation.
Let $\FmA(\PropA_1,\ldots,\PropA_k)$
 indicate that $\PropA_1,\ldots,\PropA_k$
are the propositional variables occurring
in $\FmA$ (in which case $\FmA$ is said to be $k$-ary --- unary if $k=1$, binary if $k=2$),
and let $\FmA(\FmB_1,\ldots,\FmB_k)$
refer to the formula resulting
from replacing $\FmB_i$ for each occurrence
of $\PropA_i$ in $\FmA$, for each $1 \leq i \leq k$.
Given a $\SigA$-algebra $\AlgA \SymbDef \Struct{A, \AlgInterp{\cdot}{\AlgA}}$
and a $\SigA$-formula $\FmA$,
we denote by $\AlgInterp{\FmA}{\AlgA}$
 the \emph{derived operation} 
induced on $\AlgA$ by $\FmA$.
That is, for all $\ValA_1,\ldots,\ValA_k \in A$,
provided a valuation $v$ with $v(\PropA_i) = \ValA_i$,
if $\FmA = \ConA(\FmB_1,\ldots,\FmB_k)$
we have
$\AlgInterp{\FmA}{\AlgA}(\ValA_1,\ldots,\ValA_k) = \AlgInterp{\ConA}{\AlgA}(v(\FmB_1),\ldots,v(\FmB_k))$
and, if $\FmA = \PropA_i$, then
$\AlgInterp{\FmA}{\AlgA}(\PropA_i) = v(\PropA_i)$.
By induction on the structure
of $\SigA$-formulas, we then obtain
that $\uv$ is infectious also on the derived operations of  $\AlgInfec$:

\begin{proposition}
\label{the:infectious_value_propagates}
For all $k \in \NaturalSet$, $\FmA(\PropA_1,\ldots,\PropA_k) \in \LangSetInfec$
and $\vec{\ValA} \in \ValSetInfec^k$,
$\AlgInterp{\FmA}{\AlgInfec}(\vec{\ValA})
= \AlgInterp{\FmA}{\BoolAlg}(\vec{\ValA})$
if $\vec{\ValA} \in \Set{\fv,\tv}^k$
and $\AlgInterp{\FmA}{\AlgInfec}(\vec{\ValA})
= \uv$ otherwise.
\end{proposition}

%In the present study, a \emph{logic} is a finitary \SetFmla{} or \SetSet{}
%consequence relation on $\LangSetInfec$,
%notions that we define in detail now.
In what follows, for every set $\SetA$, let $\PowerSet{\SetA}$
denote the power set of $\SetA$.
We now 
%define what are logics for us in this work.
formally introduce the notion of \emph{logic} considered in this work.

A \emph{finitary \SetSet{} consequence relation} (or a \emph{\SetSet{} logic})
over $\LangSet{\SigA}{\PropSet}$
%,
%as deeply investigated by 
%T. Shoesmith and T. Smiley
%in~\cite{ss1978},
is a binary relation $\SetSetRel$ on $\PowerSet{\LangSet{\SigA}{\PropSet}}$
satisfying 
{\ref{prop:CRSSPropO}}verlap, {\ref{prop:CRSSPropD}}ilution, {\ref{prop:CRSSPropC}}ut,
{\ref{prop:CRSSPropSS}}ubstitution-invariance
and {\ref{prop:CRSSPropF}}initariness,
for all
$\FmSetA,\FmSetB,\FmSetA',\FmSetB' \subseteq \LangSet{\SigA}{\PropSet}$:

\begin{description}[labelindent=0.5cm, labelwidth=.8cm]
    \item[\namedlabel{prop:CRSSPropO}{\CRSSPropO}] if $\FmSetA \cap \FmSetB \neq \EmptySet$, then $\FmSetA \SetSetRel{} \FmSetB$
    \item[\namedlabel{prop:CRSSPropD}{\CRSSPropD}] if $\FmSetA \SetSetRel{} \FmSetB$, then
    $\FmSetA,\FmSetA^\prime \SetSetRel{}{} \FmSetB,\FmSetB^\prime$
    \item[\namedlabel{prop:CRSSPropC}{\CRSSPropC}] if
    $\CutSetSS, \FmSetA \SetSetRel{} \FmSetB,\FmSetCompl{\CutSetSS}$
    for all $\CutSetSS\subseteq\LangSet{\SigA}{\PropSet}$,
    then $\FmSetA \SetSetRel{} \FmSetB$
    \item[\namedlabel{prop:CRSSPropSS}{\CRSSPropSS}]
    if $\FmSetA \SetSetRel{} \FmSetB$,
    then $\sigma[\FmSetA] \SetSetRel{} \sigma[\FmSetB]$,
    for every $\SigA$-substitution
    $\sigma$
    \item[\namedlabel{prop:CRSSPropF}{\CRSSPropF}] 
    if $\FmSetA \SetSetRel{} \FmSetB$,
    then 
    $\FinSet{\FmSetA} \SetSetRel{} \FinSet{\FmSetB}$
    for some finite
    $\FinSet{\FmSetA} \subseteq \FmSetA$ and
    $\FinSet{\FmSetB} \subseteq \FmSetB$
\end{description}
 \SetSet{} consequence relations have been thoroughly investigated by 
T. Shoesmith and T. Smiley
in the book~\cite{ss1978}, to which we refer the reader for further background and details.

%\noindent 

A \emph{finitary \SetFmla{} consequence
relation} (or a \emph{\SetFmla{} logic}) %, in turn, 
over $\LangSet{\SigA}{\PropSet}$
is a 
relation $\SetFmlaRel \;\subseteq\PowerSet{\LangSet{\SigA}{\PropSet}} \times \LangSet{\SigA}{\PropSet}$ satisfying the
well-known Tarskian properties
of reflexivity, monotonicity,
transitivity, substitution-invariance and finitariness.
%One may check that
\SetFmla{} logics are a particular
case of \SetSet{} logics.
One may \REV{1.4}{further} check that each \SetSet{} logic $\SetSetRel$
determines a \SetFmla{} logic 
$\SetFmlaRel_\SetSetRel{}$
over
$\LangSet{\SigA}{\PropSet}$
such that
$\FmSetA \SetFmlaRel_\SetSetRel{} \FmB$
if, and only if,
$\FmSetA \SetSetRel \Set{\FmB}$,
which is called the \emph{\SetFmla{} companion of $\SetSetRel$}.
Pairs of the form $(\FmSetA,\FmSetB)$
or $(\FmSetA,\FmB)$
are dubbed \emph{statements},
and the statements belonging
to a logic
%and assertions of the form
%$\FmSetA \SetSetRel \FmSetB$
%or
%$\FmSetA \SetSetRel \FmB$
are called \emph{consecutions} (of that logic).

A \emph{$\SigA$-matrix} is a structure $\MatA \SymbDef \Struct{\AlgA, \DSet}$,
where $\AlgA$ is a $\SigA$-algebra
and $\DSet \subseteq \ValSetA$.
We write $\overline{\DSet}$
for the set-theoretic complement $\ValSetA{\setminus}\DSet$.
The homomorphisms from $\LangAlg{\SigA}{\PropSet}$
into $\AlgA$ are called \emph{$\MatA$-valuations}.
Every $\SigA$-matrix
$\MatA$ determines a \SetSet{} consequence relation
$\SetSetRel_\MatA$
over
$\LangSet{\SigA}{\PropSet}$
%\subseteq \PowerSet{\LangSetInfec} \times \PowerSet{\LangSetInfec}$
such that
\REV{1.5}{
\begin{gather*}
\FmSetA \SetSetRel_\MatA \FmSetB
\text{ if, and only if,}
\text{ there is no $\MatA$-valuation }
\ValuationA
\text{ satisfying }\\
\ValuationA[\FmSetA] \subseteq \DSet
\text{ and }
\ValuationA[\FmSetB] \subseteq \overline\DSet.
\end{gather*}
}

We denote by $\SetFmlaRel_{\MatA}$
the \SetFmla{} companion of $\SetSetRel_{\MatA}$,
which matches the canonical
\SetFmla{} consequence relation 
over
$\LangSet{\SigA}{\PropSet}$
induced by $\MatA$,
that is,
$\FmSetA \SetFmlaRel_\MatA \FmB$
if, and only if,
there is no $\MatA$-valuation
$\ValuationA$
satisfying
$\ValuationA[\FmSetA] \subseteq \DSet$
and
$\ValuationA(\FmB) \in \overline\DSet$.
As expected, the $\SigInfec$-matrix
$\CLMat \SymbDef \Struct{\BoolAlg, \Set{\tv}}$
determines the \SetSet{} and \SetFmla{} consequence relations corresponding to
classical logic, which we  denote respectively
by $\SetSetRel_{\CLName}$ and $\SetFmlaRel_{\CLName}$.

Consider the
$\SigInfec$-matrices
 $\PWKMat \SymbDef \Struct{\AlgInfec, \Set{\uv,\tv}}$ and $\BKMat \SymbDef\Struct{\AlgInfec, \Set{\tv}}$.
Then
%, according to~\cite{pwk, bk},
Paraconsistent Weak Kleene ($\PWKName$)
and Bochvar-Kleene ($\BKName$) logics are defined, respectively,
as the \SetFmla{} logics $\SetFmlaRel_{\PWKMat}$
and $\SetFmlaRel_{\BKMat}$, which we write
$\SetFmlaRel_{\PWKName}$
and
$\SetFmlaRel_{\BKName}$ for brevity.
%In this work, 
We will also be interested
in the \SetSet{} logics
determined by these matrices
%, that is,
($\SetSetRel_{\PWKMat}$ and
$\SetSetRel_{\BKMat}$)
which we denote  simply by
$\SetSetRel_{\PWKName}$
and
$\SetSetRel_{\BKName}$,
respectively.
We may refer to them as
the
\emph{\SetSet{} versions} of $\PWKName{}$
and $\BKName{}$.

In what follows,
given a \SetFmla{} logic $\SetFmlaRel$,
we say that $\FmSetA\subseteq\LangSet{\SigA}{\PropSet}$ is \emph{$\SetFmlaRel${-}explosive} in case $\FmSetA \SetFmlaRel{} \FmA$
for all $\FmA \in \LangSet{\SigA}{\PropSet}$.
As mentioned earlier, it is well-known that
$\PWKName{}$ and
$\BKName{}$ are, respectively, the \emph{left variable inclusion companion} and the \emph{right variable inclusion companion} of classical logic, in the sense expressed by the following facts \REV{1.6}{(see~\cite{Bonzio2022book,bonzio2021plonka,caleiro2020infec} for general definitions and results concerning \emph{inclusion logics})}.

\begin{theorem}[\cite{bonzio2021plonka}]
	\label{the:pwkcltheorem}
	Let $\FmSetA, \{\FmA,\FmB\} \subseteq \LangSetInfec$.
	Then the following are equivalent:
	\begin{enumerate}
		\item $\FmSetA \SetFmlaRel_{\PWKName} \FmB$.
		\item There is $\FmSetA'\subseteq\FmSetA$
  with $\Props[\FmSetA'] \subseteq \Props(\FmB)$ such that $\FmSetA' \SetFmlaRel_{\CLName} \FmB$.
	\end{enumerate}
\end{theorem}

\begin{theorem}[\cite{urquhart2001}, Theorem 2.3.1]
	\label{the:bkcltheorem}
	Let $\FmSetA, \{\FmA,\FmB\} \subseteq \LangSetInfec$.
	Then the following are equivalent:
	\begin{enumerate}
		\item $\FmSetA \SetFmlaRel_{\BKName} \FmB$.
		\item $\FmSetA \SetFmlaRel_{\CLName} \FmB$ and $\Props(\FmB) \subseteq \Props[\FmSetA]$, or else
  $\FmSetA$ is $\SetFmlaRel_{\CLName}$-explosive.
    %$\FmSetA \SetFmlaRel_{\BKName} \FmA \land \neg\FmA$.
	\end{enumerate}
\end{theorem}

\section{Basics of Hilbert-style axiomatizations}\label{sec:Hilbertcal}
Logical matrices are a semantical way to
define %obtain concrete 
\SetFmla{} and 
\SetSet{} 
logics. %via semantics.
Another popular way are %is via 
\emph{proof systems},
which
%that is, by the manipulation of 
manipulate
syntactical objects
envisaging the construction of derivations
that bear witnesses to consecutions.
Proof systems can be classified
with respect to the \emph{proof formalism} they belong to,
based mainly on the objects they manipulate
and the shape of their rules of inference and derivations.
Each proof system induces a logic
based on the derivations one may
build via its rules of inference.

In this work, we are interested in
\emph{Hilbert-style proof systems},
or \emph{H-systems} for short.
As main characteristics,
these 
% systems 
have %, as main characteristics,
(a) rules of inference with the same shape of the 
consecutions of the induced logic;
(b) derivations as trees labelled
with sets of formulas;
and (c) the fact that they represent
a logical basis for the logics they induce, meaning that the latter is the
least logic containing the rules of
inference of the system~\cite{W88}.

Before the work of Shoesmith and
Smiley~\cite{ss1978}, rules of inference in
H-systems were constrained to
be \SetFmla{} statements, that is, pairs $(\FmSetC,\FmD) \in \PowerSet{\LangSet{\SigA}{\PropSet}} \times \LangSet{\SigA}{\PropSet}$,
usually denoted by $\HRule{\FmSetC}{\FmD}{}$,
where $\FmSetC$ is called
the \emph{antecedent} and $\FmD$,
the \emph{succedent} of the rule.
For this reason, we call them \SetFmla{} \emph{rules of inference}
and sets thereof constitute \SetFmla{} or \emph{traditional H-systems}.
They are also referred to
as \emph{single-conclusion} H-systems.
In the above-mentioned work, H-systems were
generalized to allow for multiple formulas
in the succedent of rules of inference.
In other words, rules of inference
became \SetSet{} statements, that is, pairs 
of the form
$(\FmSetC,\FmSetD) \in \PowerSet{\LangSet{\SigA}{\PropSet}} \times \PowerSet{\LangSet{\SigA}{\PropSet}}$, which we usually denote by
$\HRule{\FmSetC}{\FmSetD}{}$.
Collections of these so-called
\emph{\SetSet{} rules of inference}
constitute what we refer to as \emph{\SetSet{}} or \emph{multiple-conclusion} \emph{H-systems}.

In both formalisms, rules of inference are
usually presented \emph{schematically}, that is,
as being induced by applying 
$\SigA$-substitutions over
representative rules called \emph{rule schemas}.
An H-system is \emph{finite} when it
is presented via a finite number of rule schemas.

Users of traditional H-systems are
accustomed to derivations that are
sequences of formulas, where each
member is either a premise or results
from the application of a rule of
inference of the H-system on previous formulas in
the sequence. A proof in a traditional H-system  $\SFHSystemA$ of a 
statement $(\FmSetA,\FmB)$
is then a derivation
where the set of premises is $\FmSetA$
and the last formula is $\FmB$. Equivalently, we could
see these derivations as rooted labelled linear
trees whose nodes are labelled with sets of formulas,
where the root node is labelled with
the set of premises
and the child of each non-leaf node $\NodeA$ 
is labelled with the label $\FmSetC$ of $\NodeA$ plus the succedent of a rule of inference of $\SFHSystemA$
whose antecedent is contained in $\FmSetC$.
%any other node has label resulting from
%its parent's label by adding
%to it the formula in the succedent of a rule of inference whose antecedent is contained in that very label.
A proof of $(\FmSetA,\FmB)$,
then, is just a linear tree
whose root node is labelled with
$\FmSetA$ (or a subset thereof)
and whose leaf node contains $\FmB$.

Every \SetFmla{} H-system $\SFHSystemA$
induces a \SetFmla{} logic
$\SetFmlaRel_{\SFHSystemA}$
such that
$\FmSetA \SetFmlaRel_{\SFHSystemA} \FmB$
if and only if there is a proof of
$(\FmSetA, \FmB)$ in $\SFHSystemA$.
Given a \SetFmla{} logic
$\SetFmlaRel$
and a \SetFmla{} H-system $\SFHSystemA$,
we say that $\SFHSystemA$ is \emph{sound
for $\SetFmlaRel$} when $\SetFmlaRel_{\SFHSystemA} \, \subseteq \, \SetFmlaRel$;
that $\SFHSystemA$ is \emph{complete for $\SetFmlaRel$} when 
$\SetFmlaRel \, \subseteq \, \SetFmlaRel_{\SFHSystemA}$;
and that $\SFHSystemA$ \emph{axiomatizes $\SetFmlaRel$} (or is an \emph{axiomatization of}) $\SetFmlaRel$ when it is both sound
and complete for $\SetFmlaRel$, that is, when $\SetFmlaRel \;=\; \SetFmlaRel_{\SFHSystemA}$.

\begin{example}
\label{example:proofs-in-cl-setfmla}
The following is a well-known \SetFmla{} axiomatization
of classical logic in the signature $\SigA_{\to\neg}$, which we call $\SFHCL$ (note that it is presented by four rule schemas):

            \begin{gather*}
                \small
                \HRule{\EmptySet}{\PropB \to (\PropA \to \PropB)}{\HRuleName{\SFHCL}{1}}
                \quad
                \HRule{\EmptySet}{(\PropA\to(\PropB \to \PropC))\to((\PropA\to\PropB)\to(\PropA\to\PropC))}{\HRuleName{\SFHCL}{2}}\\[.5em]
                \HRule{\EmptySet}{(\neg\PropC\to\neg\PropB)\to((\neg\PropC\to \PropB)\to\PropC)}{\HRuleName{\SFHCL}{3}}
                \quad
                \HRule{\PropA,\PropA\to\PropB}{\PropB}{\HRuleName{\SFHCL}{4}}\\
            \end{gather*}

\noindent Here is a proof in $\SFHCL$ bearing witness
to 
%$\neg(\PropA\land\PropB) \SetFmlaRel_{\SFHCL} \neg\PropA \lor \neg\PropB$:
$\EmptySet \SetFmlaRel_{\SFHCL} \PropA \imp \PropA$:

        {
            \begin{align*}
                1. \quad & (\PropA \to ((\PropA\to\PropA) \to \PropA)) \to ((\PropA\to(\PropA\to\PropA)) \to (\PropA\to\PropA)) & \HRuleName{\SFHCL}{2}\\
                2. \quad & \PropA \to ((\PropA\to\PropA) \to \PropA) & \HRuleName{\SFHCL}{1}\\
                3. \quad & (\PropA\to(\PropA\to\PropA)) \to (\PropA\to\PropA) & 1,2,\HRuleName{\SFHCL}{4}\\
                4. \quad & \PropA\to(\PropA\to\PropA) & \HRuleName{\SFHCL}{1}\\
                5. \quad & \PropA\to\PropA & 3,4,\HRuleName{\SFHCL}{4}
            \end{align*}
            }
\end{example}

The passage from \SetFmla{} H-systems to \SetSet{} H-systems 
demands an adaptation on the latter
notions of derivations and proofs.
Now a non-leaf node $\NodeA$ may have a single child labelled with $\star$ (a \emph{discontinuation symbol}) when there is a rule of inference in the H-system
with empty succedent and whose antecedent $\FmSetC$ is contained in the label of $\NodeA$.
\REV{1.7}{This symbol indicates that the node does not
need further development (see Example~\ref{ex:derivations-classical}).}
It may alternatively have
$m$ child nodes $\NodeA_1,\ldots,\NodeA_m$ when there is
 a rule of inference $\HRule{\FmSetC}{\FmB_1,\ldots,\FmB_m}{}$ in the H-system
whose antecedent $\FmSetC$ is, as in the previous case, contained in the label of $\NodeA$.
The label of each $\NodeA_i$, in this situation, is the label of $\NodeA$
union $\Set{\FmB_i}$, for all $1 \leq i \leq m$.
See Figure~\ref{fig:derivationscheme}
for a general scheme of these derivations.
A proof of a statement $(\FmSetA,\FmSetB)$ in a \SetSet{} H-system,
then, is a labelled rooted tree
whose root node is labelled with
$\FmSetA$ (or a subset thereof)
and whose leaf nodes (now there may be more than one) 
are labelled either with $\star$
or with a set having a nonempty
intersection with $\FmSetB$.

Note that
\SetSet{} H-systems generalize
\SetFmla{} H-systems because when all rules of inference
in a \SetSet{} H-system have a single formula in the
conclusion (that is, they are \SetFmla{} rules), the derivations in that system will always be rooted labelled \emph{linear} trees, which matches our definition of \SetFmla{} derivations.

Every \SetSet{} H-system $\SSHSystemA$
induces a \SetSet{} logic
$\SetSetRel_{\SSHSystemA}$
such that
$\FmSetA \SetSetRel_{\SSHSystemA} \FmSetB$
if and only if there is a proof of
$(\FmSetA, \FmSetB)$ in $\SSHSystemA$.
Given a \SetSet{} logic
$\SetSetRel$
and a \SetSet{} H-system $\SSHSystemA$,
the notions of $\SSHSystemA$ being sound, complete or an axiomatization for $\SetSetRel$
are defined analogously as in the \SetFmla{} case.
%we say that $\SSHSystemA$ \emph{axiomatizes} $\SetSetRel$
%(or is a \SetSet{} \emph{axiomatization for}
%$\SetSetRel$) when $\SetSetRel =\SetSetRel_{\SSHSystemA}$.

\begin{figure}[t]
    \centering
    \begin{tikzpicture}[every tree node/.style={},
       level distance=1.2cm,sibling distance=1cm,
       edge from parent path={(\tikzparentnode) -- (\tikzchildnode)}, baseline]
        \Tree[.\node[style={draw,circle}] {};
            \edge[dashed];
            [.\node[style={}] {$\FmSetA$};
            \edge node[auto=right] {$\HRule{\FmSetC}{\EmptySet}{}$};
            [.{$\StarLabel$}
            ]
            ]
        ]
    \end{tikzpicture}
    \qquad
    \begin{tikzpicture}[every tree node/.style={},
       level distance=1cm,sibling distance=.5cm,
       edge from parent path={(\tikzparentnode) -- (\tikzchildnode)}, baseline]
        \Tree[.\node[style={draw,circle}] {};
            \edge[dashed];
            [.\node[style={}] {$\FmSetA$};
            \edge node[auto=right] {$\HRule{\FmSetC}{\FmD_1,\FmD_2,\ldots,\FmD_n}{}$};
            [.${\color{gray}\FmSetA},\FmD_1$
                \edge[dashed];
                [.\node[style={draw,circle}]{};
                ]
                \edge[dashed];
                [.\node[style={}]{$\cdots$};
                ]
                \edge[dashed];
                [.\node[style={draw,circle}]{};
                ]
            ]
            [.${\color{gray}\FmSetA},\FmD_2$
                \edge[dashed];
                [.\node[style={draw,circle}]{};
                ]
                \edge[dashed];
                [.\node[style={}]{$\cdots$};
                ]
                \edge[dashed];
                [.\node[style={draw,circle}]{};
                ]
            ]
            [.$\ldots$
            ]
            [.${\color{gray}\FmSetA},\FmD_n$
                \edge[dashed];
                [.\node[style={draw,circle}]{};
                ]
                \edge[dashed];
                [.\node[style={}]{$\cdots$};
                ]
                \edge[dashed];
                [.\node[style={draw,circle}]{};
                ]
            ]
            ]
        ]
    \end{tikzpicture}
    \caption{Graphical representation of $\SSHSystemA$-derivations,
    where $\SSHSystemA$ is a
    $\SetSet{}$ system.
    The dashed edges and blank circles represent other branches that may exist in the derivation.
    %in the case of expanded
    %nodes,
    We usually omit the formulas inherited from the parent
    node, 
    exhibiting only the ones introduced by the
    applied rule of inference.
    Recall that, in both cases,
    we must have $\FmSetC \subseteq \FmSetA$.}
    \label{fig:derivationscheme}
\end{figure}

\begin{example}
\label{ex:derivations-classical}
The following is a \SetSet{} axiomatization for classical logic in the signature $\SigInfec$ (in its \SetSet{} version). See Figure~\ref{fig:proofs-in-cl-setset} for examples of derivations.
\begin{gather*}
	\HRule{\EmptySet}{\PropA,\neg\PropA}{\HRuleName{\CLName}{1}} \quad
	\HRule{\PropA,\neg\PropA}{\EmptySet}{\HRuleName{\CLName}{2}} \quad
	%%%
	\HRule{\PropA,\PropB}{\PropA\land\PropB}{\HRuleName{\CLName}{3}}\quad
	\HRule{\PropA\land\PropB}{\PropA}{\HRuleName{\CLName}{4}}\quad
	\HRule{\PropA\land\PropB}{\PropB}{\HRuleName{\CLName}{5}}\\
	%%%
	\HRule{\PropA}{\PropA\lor\PropB}{\HRuleName{\CLName}{6}}\quad
	\HRule{\PropB}{\PropA\lor\PropB}{\HRuleName{\CLName}{7}}\quad
	\HRule{\PropA\lor\PropB}{\PropA,\PropB}{\HRuleName{\CLName}{8}}
\end{gather*}
%\noindent Figure~\ref{fig:proofs-in-cl-setset} shows proofs in $\SSHCL$ bearing witness, respectively,
%to 
%$\EmptySet \SetSetRel_{\SSHCL} \PropA \imp \PropA$ and
%$\neg(\PropA\land\PropB) \SetSetRel_{\SSHCL} \neg\PropA, \neg\PropB$.

        \begin{figure}[t]
        \centering
			\begin{tikzpicture}[every tree node/.style={},
						level distance=1.2cm,sibling distance=.5cm,
						edge from parent path={(\tikzparentnode) -- (\tikzchildnode)}, baseline]
			\Tree[.\node[] {$\EmptySet$};
				\edge[] node[auto=right]{${\HRuleName{\CLName}{1}}$};
				[.$\PropA$
                        \edge[] node[auto=right]{$\HRuleName{\CLName}{7}$};
                        [.$\neg\PropA\lor\PropA$
                        ]
                    ]
				[.$\neg\PropA$
                            \edge[] node[auto=right]
                            {$\HRuleName{\CLName}{6}$};
                            [.$\neg\PropA\lor\PropA$
                            ]
                    ]
			]
		\end{tikzpicture}	\qquad\qquad
            \begin{tikzpicture}[every tree node/.style={},
						level distance=1.2cm,sibling distance=.5cm,
						edge from parent path={(\tikzparentnode) -- (\tikzchildnode)}, baseline]
			\Tree[.\node[] {$\neg(\PropA\land\PropB)$};
				\edge[] node[auto=right]{$\HRuleName{\CLName}{1}$};
				[.$\PropA$
                        \edge[] node[auto=right]{$\HRuleName{\CLName}{1}$};
                        [.$\PropB$
                            \edge[] node[auto=right]{$\HRuleName{\CLName}{3}$};
                            [.$\PropA\land\PropB$
                                \edge[] node[auto=right]{$\HRuleName{\CLName}{2}$};
                                [.$\star$
                                ]
                            ]
                        ]
                        [.$\neg\PropB$
                        ]
                    ]
				[.$\neg\PropA$
                    ]
			]
		\end{tikzpicture}
        \caption{Proofs in $\SSHCL$ bearing witness, respectively,
to 
$\EmptySet \SetSetRel_{\SSHCL} \PropA \imp \PropA$ and
$\neg(\PropA\land\PropB) \SetSetRel_{\SSHCL} \neg\PropA, \neg\PropB$.}
        \label{fig:proofs-in-cl-setset}
  \end{figure}

\end{example}

The derivations 
%witnessing $\neg(\PropA\land\PropB) \SetSetRel_{\SSHCL} \neg\PropA, \neg\PropB$ 
shown in Figure~\ref{fig:proofs-in-cl-setset} have an important
property: only subformulas of
the formulas in
%$\{\neg(\PropA\land\PropB), \neg\PropA,\neg\PropB\}$ (those in 
the respective statements \REV{2.4}{$({\FmSetA},{\FmSetB})$} being proved appear in the labels of the nodes. 
%The same observation applies to
%the witness of $\EmptySet \SetSetRel_{\SSHCL} \PropA \imp \PropA$.
In fact, every statement that is provable in 
$\SSHCL$ has a proof with such feature.
For this reason, we say that $\SSHCL$ is \emph{analytic}.
Traditional (\SetFmla{}) H-systems have been historically avoided
in tasks involving proof search, as they rarely
satisfy the property of analyticity (note how the non-analyticity of $\SFHCL$ shows up in Example~\ref{example:proofs-in-cl-setfmla}). 
The solution has usually been to employ  another
deductive formalism, usually one with more meta-linguistic
resources, allowing one to prove 
meta-results that guarantee analyticity (a typical example being cut elimination in sequent-style systems~\cite{negri2001}).

%After 
Recent work by C.~Caleiro
and S.~Marcelino~\cite{marcelino19woll,marcelino19syn} %, we have 
demonstrates that the much simpler passage to \SetSet{} H-systems is enough
to obtain analytic proof systems (and thus bounded proof search) for a plethora of 
non-classical
%many-valued 
logics.
This observation will be key to us, for we will be able
to apply the %underlying 
techniques developed in the above-mentioned studies
to provide 
%the first 
finite H-systems for $\PWKName{}$ and $\BKName$.
This result, however,
 demands a slight generalization of
the notion of analyticity
in addition to the already mentioned modification of the proof formalism
to \SetSet{}.
In order to understand it,
consider first
%first let $\Theta$ be 
a set $\Theta$ of formulas on a single 
propositional variable, and let $\SSHSystemA$ be a \SetSet{} system.
The main idea
is to allow for not only subformulas of a statement %of interest 
to appear in an analytic proof,
but also formulas resulting from substitutions of those
subformulas over the formulas in $\Theta$. 
For example, if $\Theta = \{\PropC, \neg\PropC\}$,
a $\Theta$-analytic proof
witnessing that $\neg\PropA$
follows from $\neg(\PropA\land\PropB)$ would
use only formulas in $\{\PropA,\PropB,\neg\PropA,\neg\PropB,\neg\neg\PropA,\PropA\land\PropB,\neg(\PropA\land\PropB),\neg\neg(\PropA\land\PropB)\}$.
Formally,
we say that $\SSHSystemA$ is \emph{$\Theta$-analytic} whenever
$\FmSetA \SetSetRel_{\SSHSystemA} \FmSetB$ implies
that there is a \emph{$\Theta$-analytic proof} of $(\FmSetA, \FmSetB)$
in $\SSHSystemA$, that is, a proof whose nodes are labelled 
only with formulas in the set
$\Subf(\FmSetA \cup \FmSetB) \cup \Set{\FmA(\FmB) \mid \FmA \in \Theta, 
\FmB \in \Subf(\FmSetA \cup \FmSetB})$, i.e.~the \emph{$\Theta$-subformulas} of $(\FmSetA, \FmSetB)$.

%With these notions, we have that 
One can show that any finite logical matrix\footnote{Actually, the result applies to a much more general
scenario, which is not needed in the present work: the matrix can even be partial non-deterministic~\cite{marcelino19woll} in
the sense of~\cite{baaz2013,avron2011}. It may also be infinite, but then the generated system might be infinite as well.}
satisfying a very mild %inclusive
expressiveness requirement is effectively axiomatized by a 
finite
$\Theta$-analytic \SetSet{} system, for some finite $\Theta$.
This requirement is called \emph{monadicity} (or \emph{sufficient expressivess}), and intuitively means that every truth value
of the matrix can be 
%characterized in a specific sense
described by formulas
on a single variable (the set of these formulas will be 
precisely %the set 
$\Theta$).
Let us make this notion precise  and  formally state the axiomatization result.
\REV{1.8}{We say that a matrix $\MatA \SymbDef \langle \AlgA, D \rangle$ is \emph{monadic} whenever 
for every pair of \REV{2.5}{distinct} truth values $x, y \in \ValSetA$
there is a formula $\SepA$ in one propositional variable such that
$\SepA_{\AlgA}(x) \in \DSet$ and $\SepA_{\AlgA}(y) \in \ValSetA{\setminus}\DSet$
or vice-versa.} These formulas are called \emph{separators}.
Then we have that:

\begin{theorem}[\cite{marcelino19syn}, Theorem 3.5]
    \label{fact:axiomat-monadic-matrices}
    For every finite monadic logical matrix $\MatA$,
    the logic
    $\SetSetRel_{\MatA}$
    is axiomatized by a finite $\Theta$-analytic
    \SetSet{} system (which we call $\SymCalcName{\MatA}^{\Theta}$) where $\Theta$ is a finite set
    of separators for every pair of truth values of $\MatA$.
\end{theorem}

The next lemma shows why this result is so important for us in the present context.

\begin{lemma}
    \label{fact:bk-pwk-monadic}
    The matrices $\PWKMat$ and $\BKMat$ 
    are monadic, with set of separators $\Theta \SymbDef \Set{\PropA,\neg\PropA}$.
\end{lemma}
\begin{proof}
    In both matrices, $\PropA$ is a separator
    for $(\tv, \fv)$.
    In $\PWKMat$, 
    the same formula separates $(\fv,\uv)$
    and $\neg\PropA$ separates $(\uv,\tv)$.
    In $\BKMat$, 
    we have that $\PropA$ separates $(\tv,\uv)$
    and $\neg\PropA$ separates $(\fv,\uv)$.
\end{proof}

The above fact anticipates that we will be able to
provide finite and $\Set{\PropA,\neg\PropA}$-analytic \SetSet{} 
systems for the
\SetSet{} versions of $\PWKName{}$ and $\BKName{}$.
\REV{2.6}{However, it is not obvious} how to obtain
traditional finite H-systems for the original 
(and most studied) \SetFmla{} versions
of these logics. In the next couple of sections, we will
not only exhibit the announced \SetSet{} systems,
but also show how to use them to
obtain finite \SetFmla{} H-systems for $\PWKName{}$ and $\BKName{}$,
thus solving the  question regarding their finite axiomatizability.

\section{Finite H-systems for PWK}\label{sec:pwk}
Let us begin with the task of axiomatizing the \SetSet{} version of $\PWKName{}$.
The following \SetSet{} system was generated 
from the matrix $\PWKMat$
by the algorithm and simplification procedures
described in~\cite{marcelino19syn}
and implemented in~\cite[Appendix A]{greatimsc2021}, using $\Set{\PropA,\neg\PropA}$ as a set of separators (in view of Lemma~\ref{fact:bk-pwk-monadic}).
\begin{definition}
	Let $\SymCalcName{\PWKName}$ be the \SetSet{} system presented by the following rule schemas:
	\begin{gather*}
		\HRule{\EmptySet}{\PropA,\neg\PropA}{\HRuleSSName{\PWKName}{1}} \quad
		\HRule{\PropA}{\neg\neg\PropA}{\HRuleSSName{\PWKName}{2}} \quad
		\HRule{\neg\neg\PropA}{\PropA}{\HRuleSSName{\PWKName}{3}}\\
		%%%
		\HRule{\PropA,\PropB}{\PropA\land\PropB}{\HRuleSSName{\PWKName}{4}}\quad
		\HRule{\PropA\land\PropB}{\PropA,\PropB}{\HRuleSSName{\PWKName}{5}}\quad
		\HRule{\PropA\land\PropB}{\PropA,\neg\PropB}{\HRuleSSName{\PWKName}{6}}\quad
		\HRule{\PropA\land\PropB}{\neg\PropA,\PropB}{\HRuleSSName{\PWKName}{7}}\quad
		\HRule{\neg(\PropA\land\PropB)}{\neg\PropA,\neg\PropB}{\HRuleSSName{\PWKName}{8}}\\
		\HRule{\PropA,\neg\PropA}{\PropA\land\PropB}{\HRuleSSName{\PWKName}{9}}
        \quad
		%\HRule{\PropA,\neg\PropA}{\neg(\PropA\land\PropB)}{\HRuleSSName{\PWKName}{10}}
		\REV{3.1}{\HRule{\neg\PropA}{\neg(\PropA\land\PropB)}{\HRuleSSName{\PWKName}{10}}}
        \quad
		\HRule{\PropB,\neg\PropB}{\PropA\land\PropB}{\HRuleSSName{\PWKName}{11}}
        \quad
		%\HRule{\PropB,\neg\PropB}{\neg(\PropA\land\PropB)}{\HRuleSSName{\PWKName}{12}}
		\REV{3.1}{\HRule{\neg\PropB}{\neg(\PropA\land\PropB)}{\HRuleSSName{\PWKName}{12}}}
        \\
		%%%
		\HRule{\PropA}{\PropA\lor\PropB}{\HRuleSSName{\PWKName}{13}}\quad
		\HRule{\PropB}{\PropA\lor\PropB}{\HRuleSSName{\PWKName}{14}}\quad
		\HRule{\PropA\lor\PropB}{\PropA,\PropB}{\HRuleSSName{\PWKName}{15}}\quad
		\HRule{\neg(\PropA\lor\PropB)}{\PropA,\neg\PropB}{\HRuleSSName{\PWKName}{16}}\\
		\HRule{\neg(\PropA\lor\PropB)}{\neg\PropA,\PropB}{\HRuleSSName{\PWKName}{17}}\quad
		\HRule{\neg(\PropA\lor\PropB)}{\neg\PropA,\neg\PropB}{\HRuleSSName{\PWKName}{18}}\quad
		\HRule{\PropA,\neg\PropA}{\neg(\PropA\lor\PropB)}{\HRuleSSName{\PWKName}{19}}\quad
		\HRule{\PropB,\neg\PropB}{\neg(\PropA\lor\PropB)}{\HRuleSSName{\PWKName}{20}}
	\end{gather*}

% {\color{red} MIssing
%$$\frac{}{p,\neg(p\land q)}\qquad \frac{}{q,\neg(p\land q)}$$
%or enough the refinement of the ones we have
%$$\frac{\neg p}{\neg(p\land q)}\qquad \frac{\neg q}{\neg(p\land q)}$$
%
%and
%
%$$\frac{\neg(p\lor q)}{\neg p}\qquad \frac{\neg(p\lor q)}{\neg q}$$
%
%
%}
\end{definition}

Since this system is equivalent to
the system $\SymCalcName{\PWKMat}^{\Set{\PropA,\neg\PropA}}$
mentioned in Theorem~\ref{fact:axiomat-monadic-matrices} (when specialized to $\PWKMat$),
and since the mentioned simplification procedures preserve $\Theta$-analyticity, we obtain:

\begin{theorem}
    $\SymCalcName{\PWKName}$ is $\Set{\PropA,\neg\PropA}$-analytic
    and
$\SetSetRel_{\SymCalcName{\PWKName}} = \SetSetRel_{\PWKName}$.
\end{theorem}
%\begin{proof}
%The matrix
%    $\PWKMat$ is monadic and
%    $\SymCalcName{\PWKName}$
%     was generated by the axiomatization procedure described
%     in~\cite{marcelino}.
%\end{proof}

Our goal now is to find a finite \SetFmla{} H-system for $\PWKName$.
We will see that this task is
easily solved because the disjunction connective in this logic allows us
to convert $\SymCalcName{\PWKName}$ into 
the desired finite \SetFmla{} system. 
More generally, 
every \SetFmla{} logic $\SetFmlaRel$ is finitely axiomatized
by a \SetFmla{} H-system whenever it satisfies two
conditions which we now describe~\cite[Theorem 5.37]{ss1978}. First, the logic is the \SetFmla{} companion
of a \SetSet{} logic finitely axiomatized by a \SetSet{} H-system, say $\SSHSystemA$.
Second, it satisfies the following property for some 
%definable 
binary formula $C(\PropA,\PropB)$ (said to be a \emph{definable binary connective} in this context):
\begin{gather*}\tag{{\normalfont\textsf{disj}}}
    \label{eq:disj-prop}
	\text{for all } \FmSetA \cup \{\FmA,\FmB,\FmC\} \subseteq \LangSet{\SigA}{\PropSet},\\
	\FmSetA, \FmA \vdash_{} \FmC
	\text{ and }
	\FmSetA, \FmB \vdash_{} \FmC
	\text{ if, and only if, }
	\FmSetA, C(\FmA,\FmB) \vdash_{} \FmC.
\end{gather*}
The proof of this fact in~\cite{ss1978} reveals how to effectively
convert $\SSHSystemA$ into the desired \SetFmla{} H-system.
Let us see how to perform this conversion and then apply the 
transformation to $\SymCalcName{\PWKName}$\REV{2.7}{.}\footnote{Note that we use $\lor$ to simplify
notation, but the same definition could be rephrased
with the derived connective $C(\PropA,\PropB)$.}
In what follows,
when $\FmSetA \SymbDef \{\FmA_1,\ldots,\FmA_n\} \subseteq \LangSetInfec$ ($n \geq 1$),
let $\bigvee\FmSetA \SymbDef (\ldots(\FmA_1 \lor \FmA_2) \lor \ldots) \lor \FmA_n$. Also, let $\FmSetA \lor \FmB \SymbDef \{\FmA \lor \FmB \mid \FmA \in \FmSetA\}$. Note that the latter set is empty when $\FmSetA$ is empty.

\begin{definition}
Let $\SymCalcName{}$ be
a \SetSet{} system
and
$p_0$ be a
    propositional variable
    not occurring in the
    rule schemas of $\mathsf{R}$.
Define the system $\mathsf{R}^{\lor}$ as being presented by the rule schemas
$\left\{\frac{p \lor p}{p},\frac{p}{p \lor q}, \frac{p \lor q}{q \lor p}, \frac{p \lor (q \lor r)}{(p \lor q) \lor r}\right\} \cup
\left\{\mathsf{r}^\lor \mid \mathsf{r} \text{ is a rule schema of } \mathsf{R}\right\}$ where $\mathsf{r}^\lor$
is $\frac{\varnothing}{\FmA}$
    if $\mathsf{r} = \frac{\varnothing}{\FmA}$, $\frac{\FmSetA \lor p_0}{(\bigvee \FmSetB) \lor p_0}$
    if $\mathsf{r} = \frac{\FmSetA}{\FmSetB}$,
    and $\frac{\FmSetA \lor p_0}{p_0}$
    if $\mathsf{r} = \frac{\FmSetA}{\varnothing}$.
\end{definition}

Below we present the result of this procedure when applied to $\SymCalcName{\PWKName}$.
Note that the conversion of rule $\HRuleSSName{\PWKName}{15}$
results in a rule of the form $\FmA/\FmA$, and thus can be discarded.

\begin{definition}
    Let $\ASymCalcName{\PWKName}$ be the \SetFmla{} system presented by the following rule schemas:

    \begin{gather*}
        \HRule{\EmptySet}{\PropA\lor\neg\PropA}{\HRuleName{\PWKName}{1}} \quad
        \HRule{\PropA\lor\PropC}{\neg\neg\PropA\lor\PropC}{\HRuleName{\PWKName}{2}} \quad
        \HRule{\neg\neg\PropA \lor \PropC}{\PropA\lor\PropC}{\HRuleName{\PWKName}{3}}\\
        %%%
        \HRule{\PropA \lor \PropC,\PropB \lor \PropC}{(\PropA\land\PropB) \lor \PropC}{\HRuleName{\PWKName}{4}}\quad
        \HRule{(\PropA\land\PropB)\lor\PropC}{(\PropA\lor\PropB)\lor\PropC}{\HRuleName{\PWKName}{5}}\quad
        \HRule{(\PropA\land\PropB)\lor\PropC}{(\PropA \lor \neg\PropB)\lor\PropC}{\HRuleName{\PWKName}{6}}\\
        %%%%
        \HRule{(\PropA\land\PropB)\lor\PropC}{(\neg\PropA\lor\PropB)\lor\PropC}{\HRuleName{\PWKName}{7}}\quad
        \HRule{\neg(\PropA\land\PropB)\lor\PropC}{(\neg\PropA \lor \neg\PropB)\lor\PropC}{\HRuleName{\PWKName}{8}}\quad
        \HRule{\PropA \lor \PropC,\neg\PropA \lor \PropC}{(\PropA\land\PropB)\lor\PropC}{\HRuleName{\PWKName}{9}}\\
        %%%%
        %\HRule{\PropA,\neg\PropA}{\neg(\PropA\land\PropB)}{\HRuleName{\PWKName}{10}}
        \REV{3.1}{\HRule{\neg\PropA \lor \PropC}{\neg(\PropA\land\PropB) \lor \PropC}{\HRuleName{\PWKName}{10}}}
        \quad
        \HRule{\PropB \lor \PropC,\neg\PropB \lor \PropC}{(\PropA\land\PropB)\lor\PropC}{\HRuleName{\PWKName}{11}}
        \quad
        %\HRule{\PropB,\neg\PropB}{\neg(\PropA\land\PropB)}{\HRuleName{\PWKName}{12}}
        \REV{3.1}{\HRule{\neg\PropB \lor \PropC}{\neg(\PropA\land\PropB) \lor \PropC}{\HRuleName{\PWKName}{12}}}
        \\
        %%%
        \HRule{\PropA \lor \PropC}{(\PropA\lor\PropB)\lor\PropC}{\HRuleName{\PWKName}{13}}\quad
        \HRule{\PropB \lor \PropC}{(\PropA\lor\PropB) \lor\PropC}{\HRuleName{\PWKName}{14}}\quad
        %\HRule{\PropC\lor(\PropA\lor\PropB)}{\PropC\lor(\PropA \lor \PropB)}{\HRuleName{\PWKName}{15}}\quad
        \HRule{\neg(\PropA\lor\PropB)\lor\PropC}{(\PropA\lor\neg\PropB)\lor\PropC}{\HRuleName{\PWKName}{15}}\\
        \HRule{\neg(\PropA\lor\PropB)\lor\PropC}{(\neg\PropA\lor\PropB)\lor\PropC}{\HRuleName{\PWKName}{16}}\quad
        \HRule{\neg(\PropA\lor\PropB)\lor\PropC}{(\neg\PropA \lor \neg\PropB)\lor\PropC}{\HRuleName{\PWKName}{17}}\\
        \HRule{\PropA\lor\PropC,\neg\PropA\lor\PropC}{\neg(\PropA\lor\PropB)\lor\PropC}{\HRuleName{\PWKName}{18}}
        \quad
        \HRule{\PropB \lor \PropC,\neg\PropB \lor \PropC}{\neg(\PropA\lor\PropB) \lor \PropC}{\HRuleName{\PWKName}{19}}\\
        \HRule{\PropA\lor\PropA}{\PropA}{\HRuleName{\PWKName}{20}}
        \quad
        \HRule{\PropA}{\PropA\lor\PropB}{\HRuleName{\PWKName}{21}}
        \quad
        \HRule{\PropA \lor \PropB}{\PropB\lor\PropA}{\HRuleName{\PWKName}{22}}\quad
        \HRule{\PropA \lor (\PropB \lor \PropC)}{(\PropA\lor\PropB)\lor\PropC}{\HRuleName{\PWKName}{23}} 
    \end{gather*}         
\end{definition}

As anticipated in the previous discussion, we have that:

\begin{theorem}[\cite{ss1978}, Theorem 5.37]
    \label{fact:setfmla-from-setset}
    If $\SetSetRel_{\SymCalcName{}} = \SetSetRel$
    and $\SetFmlaRel_\SetSetRel$ satisfies {\normalfont (\ref{eq:disj-prop})}, then
    $\SetFmlaRel_{\SymCalcName{}^\lor} \; = \; \SetFmlaRel_\SetSetRel.$
\end{theorem}

\begin{remark}
    The authors of~\cite{ss1978}
    also show that a similar conversion
    between \SetSet{} and \SetFmla{}
    is possible when the logic
    has a definable binary connective $C(\PropA,\PropB)$ that satisfies 
    the so-called \emph{deduction theorem}:
    \begin{gather*}\tag{{\normalfont\textsf{ded}}}
    \label{eq:ded-prop}
	\text{for all } \FmSetA \cup \{\FmA,\FmB\} \subseteq \LangSet{\SigA}{\PropSet},\\
	\FmSetA, \FmA \vdash_{} \FmB
	\text{ if, and only if, }
	\FmSetA \vdash_{} C(\FmA,\FmB).
\end{gather*}
\end{remark}

Theorem~\ref{fact:setfmla-from-setset}  can then be applied to $\PWKName$ because
$\SetFmlaRel_{\PWKName} \;=\; \SetFmlaRel_{\SetSetRel_{\PWKName}}$ and
it satisfies (\ref{eq:disj-prop}),
as we establish below.

\begin{proposition}
        \label{fact:pwk-satisf-disj}
	For all $\FmSetA \cup \{\FmA,\FmB,\FmC\} \subseteq \LangSetInfec$,
	\[
	\FmSetA, \FmA \vdash_{\PWKName} \FmC
	\text{ and }
	\FmSetA, \FmB \vdash_{\PWKName} \FmC
	\text{ if, and only if, }
	\FmSetA, \FmA\lor\FmB \vdash_{\PWKName} \FmC.\]
\end{proposition}
\begin{proof}
        The reader can easily check that
	 the presence of rules 
	$\HRuleSSName{\PWKName}{13}$,
	$\HRuleSSName{\PWKName}{14}$ and
	$\HRuleSSName{\PWKName}{15}$
        in $\SymCalcName{\PWKName}$
        is enough to prove this statement.
\end{proof}

%We finally reach the main result of this section
%as a direct application of 
%Theorem~\ref{fact:setfmla-from-setset} and
%Proposition~\ref{fact:pwk-satisf-disj}:

In other words,

\begin{theorem}
    $\SetFmlaRel_{\ASymCalcName{\PWKName}} \; = \; \SetFmlaRel_{\PWKName}$.
\end{theorem}

\section{Finite H-systems for BK}\label{sec:bk}
%\begin{theorem}
%	\label{the:infectious_value_propagates}
%	For all $\FmA(\PropA_1,\ldots,\PropA_k) \in L$,
%	$\FmA_{\BKMat}(\ValA_1,\ldots,\ValA_k) = \uv$
%	if, and only if, $\ValA_i = \uv$ for some $1 \leq i \leq k$.
%\end{theorem}
%\begin{proof}
%	\TODO{Mention where this was proved.}
%\end{proof}

We shall proceed as in the previous case, 
%of $\PWKName$ in the previous section by 
starting 
with the axiomatization of the \SetSet{} version
of $\BKName$.
In view of Lemma~\ref{fact:bk-pwk-monadic}, we can apply the same reasoning
as the one applied to axiomatize $\PWKName$ in \SetSet{}, that is, we can automatically generate
a $\Set{\PropA,\neg\PropA}$-analytic axiomatization for $\BKName$:

\begin{definition}
	Let $\SymCalcName{\BKName}$ be the \SetSet{} system presented by the following rule schemas:
\begin{gather*}
	\HRule{\PropA,\neg\PropA}{\EmptySet}{\HRuleSSName{\BKName}{1}} \quad
	\HRule{\PropA}{\neg\neg\PropA}{\HRuleSSName{\BKName}{2}} \quad
	\HRule{\neg\neg\PropA}{\PropA}{\HRuleSSName{\BKName}{3}}\\
	%%%
	\HRule{\PropA,\PropB}{\PropA\land\PropB}{\HRuleSSName{\BKName}{4}}\quad
	\HRule{\neg\PropA,\neg\PropB}{\neg(\PropA\land\PropB)}{\HRuleSSName{\BKName}{5}}\quad
	\HRule{\neg\PropA,\PropB}{\neg(\PropA\land\PropB)}{\HRuleSSName{\BKName}{6}}\quad
	\HRule{\PropA,\neg\PropB}{\neg(\PropA\land\PropB)}{\HRuleSSName{\BKName}{7}}\\
	\HRule{\neg(\PropA\land\PropB)}{\neg\PropA,\PropA}{\HRuleSSName{\BKName}{8}}\quad
	\HRule{\neg(\PropA\land\PropB)}{\neg\PropB,\PropB}{\HRuleSSName{\BKName}{9}}\quad
	\HRule{\PropA\land\PropB}{\PropA}{\HRuleSSName{\BKName}{10}}\quad
	\HRule{\PropA\land\PropB}{\PropB}{\HRuleSSName{\BKName}{11}}\\
	%%
	%\HRule{\neg(\PropA\lor\PropB),\PropA}{\EmptySet}{\HRuleName{\BKName}{12}}\quad
	%\HRule{\neg(\PropA\lor\PropB),\PropB}{\EmptySet}{\HRuleName{\BKName}{13}}\quad
	\HRule{\neg\PropA,\neg\PropB}{\neg(\PropA\lor\PropB)}{\HRuleSSName{\BKName}{12}}\quad
	\HRule{\neg(\PropA\lor\PropB)}{\neg\PropA}{\HRuleSSName{\BKName}{13}}\quad
	\HRule{\neg(\PropA\lor\PropB)}{\neg\PropB}{\HRuleSSName{\BKName}{14}}\quad
	\HRule{\PropA\lor\PropB}{\PropA,\neg\PropA}{\HRuleSSName{\BKName}{15}}\\
	\HRule{\PropA\lor\PropB}{\PropB,\neg\PropB}{\HRuleSSName{\BKName}{16}}\quad
	\HRule{\neg\PropA,\PropB}{\PropA\lor\PropB}{\HRuleSSName{\BKName}{17}}\quad
	\HRule{\PropA,\neg\PropB}{\PropA\lor\PropB}{\HRuleSSName{\BKName}{18}}\quad
	\HRule{\PropA,\PropB}{\PropA\lor\PropB}{\HRuleSSName{\BKName}{19}}\quad
	\HRule{\PropA\lor\PropB}{\PropA,\PropB}{\HRuleSSName{\BKName}{20}}\quad
\end{gather*}

\end{definition}

As in the case of $\PWKName$, 
since $\SymCalcName{\BKName}$ is equivalent to
the system $\SymCalcName{\BKMat}^{\Set{\PropA,\neg\PropA}}$
mentioned in Theorem~\ref{fact:axiomat-monadic-matrices} (when particularized to $\BKMat$),
and the employed simplification procedures 
preserve $\Theta$-analyticity, we have:

\begin{theorem}
    $\SymCalcName{\BKName}$ is $\Set{\PropA,\neg\PropA}$-analytic
    and
$\SetSetRel_{\SymCalcName{\BKName}} = \SetSetRel_{\BKName}$.
\end{theorem}

\begin{remark}
    It is not hard to see that $\BKMat$ 
    results from a renaming of
    the truth-values of
    %is 
    %isomorphic (for a precise definition, see~\cite{somewhere}) to
    the logical matrix
    $\MatA' \SymbDef \Struct{\AlgInfec', \Set{\fv}}$, where
    $\AlgInfec'$ has the same set $\ValSetInfec$ of truth values and
    its truth tables are such that
    $\AlgInterp{\lor}{\AlgInfec'} \SymbDef \AlgInterp{\land}{\AlgInfec}$,
    $\AlgInterp{\land}{\AlgInfec'} \SymbDef \AlgInterp{\lor}{\AlgInfec}$
    and
    $\AlgInterp{\neg}{\AlgInfec'} \SymbDef \AlgInterp{\neg}{\AlgInfec}$
    (just swap %replace 
    $\tv$ and $\fv$ %and vice-versa 
    in the interpretations and in the designated set).
    %In other words, $\BKMat$
    %is isomorphic to $\MatA'$.
    Note also that, if we take $\PWKMat$ and replace
    its designated set $\Set{\tv,\uv}$ 
    %by 
    %its complement with respect to
    %$\ValSetInfec$, that is,
    by $\Set{\fv}$
    and swap 
    %change
    the truth tables of $\land$ and $\lor$, %one by the other, 
    we obtain $\MatA'$.
    %, we 
    %get a matrix that is isomorphic to $\MatA'$ (the only change is that
    %the truth tables of $\land$ and $\lor$ were changed one by the other).
    %We have that
    The axiomatization procedure of~\cite{marcelino19syn}
    implies in this situation that $\MatA'$ is axiomatized simply by
    taking $\SymCalcName{\PWKName}$ and turning
    its rules of inference upside down (antecedents become succedents, and vice-versa),
    in addition to replacing $\land$ by $\lor$ and
    vice-versa in the rules. 
    We call the resulting system the \emph{dualization} of $\SymCalcName{\PWKName}$.
    \REV{1.10}{Because $\BKMat$ results from $\MatA'$ 
    by this simple renaming
    of truth values, we have that it is axiomatized by this same \SetSet{} system.}
    The reader can easily check that, indeed, $\SymCalcName{\BKName}$
    is just the dualization of $\SymCalcName{\PWKName}$.
\end{remark}

Finding a finite \SetFmla{} axiomatization for $\BKName$ turns out to be harder than in
the case of $\PWKName$. The reason, as we prove in the next proposition, is that in $\BKName$ it is impossible to define a 
binary connective $C(\PropA,\PropB)$
satisfying (\ref{eq:disj-prop}) 
or (\ref{eq:ded-prop}).

\begin{proposition}
        \label{fact:no-disj-no-ded-theo}
	The following holds for $\BKName$:
	\begin{enumerate}
		\item For no binary formula $C(\PropA,\PropB) \in \LangSetInfec$ we have
			$\FmSetA,\FmA \SetFmlaRel_{\BKName} \FmC$ and $\FmSetA,\FmB \SetFmlaRel_{\BKName} \FmC$
			whenever $\FmSetA, C(\FmA,\FmB) \SetFmlaRel_{\BKName} \FmC$,
			for all $\FmSetA \cup \Set{\FmA,\FmB, \FmC} \subseteq \LangSetInfec$.
		\item For no binary formula $C(\PropA,\PropB) \in \LangSetInfec$ we have
			$\FmSetA \SetFmlaRel_{\BKName} C(\FmA,\FmB)$ whenever
			$\FmSetA,\FmA \SetFmlaRel_{\BKName} \FmB$,
			for all $\FmSetA \cup \Set{\FmA,\FmB} \subseteq \LangSetInfec$.
	\end{enumerate}
\end{proposition}
\begin{proof}
	For item 1, note that $C(\PropA,\PropB) \SetFmlaRel_{\BKName} C(\PropA,\PropB) \lor \neg C(\PropA,\PropB)$,
	however  $\PropB \not\SetFmlaRel_{\BKName} C(\PropA,\PropB) \lor \neg C(\PropA,\PropB)$,
	as a $\BKName$-valuation assigning $\uv$ to $\PropA$
	and $\tv$ to $\PropB$ would be a countermodel for the latter consecution (see Theorem~\ref{the:infectious_value_propagates}).
	Similarly, for item 2,
	note that $\neg\PropA,\PropA \SetFmlaRel_{\BKName} \PropB$,
	but $\neg\PropA \not\SetFmlaRel_{\BKName} C(\PropA, \PropB)$,
	what can be seen by considering a $\BKName$-valuation
	assigning $\fv$ to $\PropA$
	and $\uv$ to $\PropB$.
\end{proof}

Therefore, up to this point, the mere existence of a finite \SetSet{} system for
$\BKName$ does not guarantee that this logic is finitely axiomatizable in \SetFmla{}.
It does not mean, however, that such system cannot help us
in an \emph{ad hoc} effort to finitely axiomatize $\BKName$.

We begin by noting that only the rules $\HRuleSSName{\BKName}{i}$, with $i \in \Set{8,9,15,16,20}$,
have multiple formulas in the succedent. We will replace the first four
of these by the following \SetFmla{} rules:
\begin{gather*}
	\HRule{\neg(\PropA\land\PropB)}{\neg\PropA\lor\PropA}{\HRuleName{\BKName{}}{8\star}}\quad
	\HRule{\neg(\PropA\land\PropB)}{\neg\PropB\lor\PropB}{\HRuleName{\BKName}{9\star}}\quad
	\HRule{\PropA\lor\PropB}{\PropA\lor\neg\PropA}{\HRuleName{\BKName}{15\star}}\quad
	\HRule{\PropA\lor\PropB}{\PropB\lor\neg\PropB}{\HRuleName{\BKName}{16\star}}\quad
\end{gather*}

\begin{definition}
	Let $\SymCalcName{\BKName\star}$ be 
	$\SymCalcName{\BKName}$ 
	but with 
	$\HRuleSSName{\BKName}{i}$
	replaced by $\HRuleName{\BKName}{i\star}$, for each $i \in \Set{8,9,15,16}$.
\end{definition}

Then we have that this transformation preserves
the induced logic:

\begin{proposition}
	$\SymCalcName{\BKName}$
	and 
	$\SymCalcName{\BKName\star}$
	induce the same \SetSet{} logic.
\end{proposition}
\begin{proof}
    We just need to show that
    $\SetSetRel_{\SymCalcName{\BKName}} = \SetSetRel_{\SymCalcName{\BKName\star}}$.
	The right-to-left inclusion is easy,
	and the converse follows thanks to the presence of $\HRuleSSName{\BKName}{20}$.
\end{proof}

%The next example shows a derivation of one of the De Morgan rules 
%in $\ASymCalcName{\BKName\star}$.

\begin{example}
	\label{ex:demorgan}
	The following derivation bears witness to
	$\neg(\PropA\land\PropB) \SetSetRel_{\SymCalcName{\BKName\star}} \neg\PropA \lor \neg\PropB$:

    \begin{figure}[H]
        \centering
			\begin{tikzpicture}[every tree node/.style={},
						level distance=1.2cm,sibling distance=.5cm,
						edge from parent path={(\tikzparentnode) -- (\tikzchildnode)}, baseline]
			\Tree[.\node[] {$\neg(\PropA\land\PropB)$};
				\edge[] node[auto=right]{$\HRuleName{\BKName}{9\star}$};
				[.$\PropB\lor\neg\PropB$
					\edge[] node[auto=right]{$\HRuleSSName{\BKName}{20}$};
					[.\node[] {$\PropB$};
						\edge[] node[auto=right]{$\HRuleSSName{\BKName}{2}$};
						[.\node[] {$\neg\neg\PropB$};
							\edge[] node[auto=right]{$\HRuleName{\BKName}{8\star}$};
							[.\node[] {$\PropA\lor\neg\PropA$};
								\edge[] node[auto=right]{$\HRuleSSName{\BKName}{20}$};
								[.$\PropA$
									\edge[] node[auto=right]{$\HRuleSSName{\BKName}{4}$};
									[.$\PropA\land\PropB$
										\edge[] node[auto=right]{$\HRuleSSName{\BKName}{1}$};
										[.$\star$
										]
									]
								]
								[.$\neg\PropA$
									\edge[] node[auto=right]{$\HRuleSSName{\BKName}{18}$};
									[.$\neg\PropA\lor\neg\PropB$
									]
								]
							]
						]
					]
					[.\node[] {$\neg\PropB$};
						\edge[] node[auto=right]{$\HRuleName{\BKName}{8\star}$};
						[.$\PropA\lor\neg\PropA$
							\edge[] node[auto=right]{$\HRuleSSName{\BKName}{20}$};
							[.$\PropA$
								\edge[] node[auto=right]{$\HRuleSSName{\BKName}{2}$};
								[.$\neg\neg\PropA$
									\edge[] node[auto=right]{$\HRuleSSName{\BKName}{17}$};
									[.$\neg\PropA\lor\neg\PropB$
									]
								]
							]
							[.$\neg\PropA$
								\edge[] node[auto=right]{$\HRuleSSName{\BKName}{19}$};
								[.$\neg\PropA\lor\neg\PropB$
								]
							]
						]
					]
				]
			]
		\end{tikzpicture}	
            \label{fig:deriv-in-bk-star}
            \caption{A derivation in $\SymCalcName{\BKName\star}$ showing that
	$\neg(\PropA\land\PropB) \SetSetRel_{\SymCalcName{\BKName\star}} \neg\PropA \lor \neg\PropB$.}
    \end{figure}
\end{example}

\begin{remark}
    The modifications in $\SymCalcName{\BKName}$ 
    that resulted in $\SymCalcName{\BKName\star}$, despite preserving the
    induced logic, are not guaranteed to preserve $\Set{\PropA,\neg\PropA}$-analyticity.
    The previous example may be seen as an illustration of this fact.
\end{remark}

The fact that the only rule of $\SymCalcName{\BKName\star}$
with more than one formula in the succedent is $\HRule{\PropA\lor\PropB}{\PropA,\PropB}{\HRuleSSName{\BKName}{20}}$
will help us in providing a finite \SetFmla{} system for $\BKName$, thus answering
positively the question of its finite axiomatizability.
Before showing why and how, let us introduce 
some transformations over \SetFmla{} rules
that will be useful in our endeavour:

\begin{definition}
	\label{def:lifted}
	Let $\HRule{\FmA_1,\ldots,\FmA_m}{\FmB}{\RuleA}$ be a \SetFmla{} inference rule
	and $\PropC$ be a propositional variable not occurring in \REV{2.8}{any} of the formulas $\FmA_1,\ldots,\FmA_m$ and $\FmB$.
 For simplicity, we  define the binary connective $\imp$
by abbreviation: for all $\FmA,\FmB \in \LangSetInfec$,
let $\FmA\imp\FmB \SymbDef \neg\FmA\lor\FmB$.
	Then:
	\begin{enumerate}
		\item The \emph{$\lor$-lifted version} of $\RuleA$
			is the rule
			\[
				\HRule{\PropC \lor \FmA_1, \ldots, \PropC \lor \FmA_m}{\PropC \lor \FmB}{\RuleA^\lor}
			\]
		\item The \emph{$\to$-lifted version} of $\RuleA$
			is the rule
			\[
				\HRule{\PropC \to \FmA_1, \ldots, \PropC \to \FmA_m}{\PropC \to \FmB}{\RuleA^\to}
			\]
	\end{enumerate}
\end{definition}

The following characterization of rules of inference will also
be useful to us, in view of Theorem~\ref{the:bkcltheorem}:

\begin{definition}
	\label{def:containment_condition}
	A \SetFmla{} inference rule $\HRule{\FmSetA}{\FmB}{\RuleA}$ 
	is said to satisfy the \emph{containment condition} whenever
	$\Props(\FmB) \subseteq \Props[\FmSetA]$.
\end{definition}

 We will provide
a \SetFmla{} H-system resulting from $\SymCalcName{\BKName\star}$ essentially by the following modifications: 
\begin{enumerate}[label=(\alph*),leftmargin=1cm]
    \item removing the rule 
$\HRuleSSName{\BKName}{20}$; 
\item replacing
$\HRuleSSName{\BKName}{1}$, a rule with empty succedent, with a new rule called $\HRuleName{\BKName}{1\star}$ having a fresh variable in the succedent;
\item adding some rules concerning $\lor$;
\item adding all $\lor$-lifted versions (see Definition~\ref{def:lifted}) of
 all rules but $\HRuleName{\BKName}{1\star}$.
\end{enumerate}
%(a) removing the rule 
%$\HRuleSSName{\BKName}{20}$; 
%(b) replacing
%$\HRuleName{\BKName}{1}$, a rule with empty succedent, with a new rule called $\HRuleName{\BKName}{1\star}$ having a %fresh variable in the succedent;
%and removing the others, as they become derivable in this new setting (see Proposition~\ref{prop:basic_deriv_rules_bk});
%(c) adding some rules concerning $\lor$;
%and (d) adding all $\lor$-lifted versions (see Definition~\ref{def:lifted}) of
%the all rules but $\HRuleName{\BKName}{1\star}$.
Having the lifted rules for all rules satisfying the
containment condition will be important for completeness, as we will see.
Our task, then, boils down to showing that applications of 
$\HRuleSSName{\BKName}{1}$
%$\HRuleName{\BKName}{12}$,
%$\HRuleName{\BKName}{13}$ 
and
$\HRuleSSName{\BKName}{20}$ in derivations in $\SymCalcName{\BKName\star}$
of \SetFmla{} statements
may be replaced by applications of rules of the proposed \SetFmla{} system.
We display this system below for clarity and ease of reference.

\begin{definition}
	\label{def:bksetfmla}
	Let $\ASymCalcName{\BKName}$ be the \SetFmla{} system given by the rule schemas
	\begin{gather*}
	\HRule{\PropA,\neg\PropA}{\PropB}{\HRuleName{\BKName}{1\star}} \quad
	\HRule{\PropA}{\neg\neg\PropA}{\HRuleName{\BKName}{2}} \quad
	\HRule{\neg\neg\PropA}{\PropA}{\HRuleName{\BKName}{3}}\\
	%%%
	\HRule{\PropA,\PropB}{\PropA\land\PropB}{\HRuleName{\BKName}{4}}\quad
	\HRule{\neg\PropA,\neg\PropB}{\neg(\PropA\land\PropB)}{\HRuleName{\BKName}{5}}\quad
	\HRule{\neg\PropA,\PropB}{\neg(\PropA\land\PropB)}{\HRuleName{\BKName}{6}}\quad
	\HRule{\PropA,\neg\PropB}{\neg(\PropA\land\PropB)}{\HRuleName{\BKName}{7}}\\
	\HRule{\neg(\PropA\land\PropB)}{\neg\PropA\lor\PropA}{\HRuleName{\BKName{}}{8\star}}\quad
	\HRule{\neg(\PropA\land\PropB)}{\neg\PropB\lor\PropB}{\HRuleName{\BKName}{9\star}}\quad
	\HRule{\PropA\land\PropB}{\PropA}{\HRuleName{\BKName}{10}}\quad
	\HRule{\PropA\land\PropB}{\PropB}{\HRuleName{\BKName}{11}}\\
	\HRule{\neg\PropA,\neg\PropB}{\neg(\PropA\lor\PropB)}{\HRuleName{\BKName}{12}}\quad
	\HRule{\neg(\PropA\lor\PropB)}{\neg\PropA}{\HRuleName{\BKName}{13}}\quad
	\HRule{\neg(\PropA\lor\PropB)}{\neg\PropB}{\HRuleName{\BKName}{14}}\quad
	\HRule{\PropA\lor\PropB}{\PropA\lor\neg\PropA}{\HRuleName{\BKName}{15\star}}\\
	\HRule{\PropA\lor\PropB}{\PropB\lor\neg\PropB}{\HRuleName{\BKName}{16\star}}\quad
	\HRule{\neg\PropA,\PropB}{\PropA\lor\PropB}{\HRuleName{\BKName}{17}}\quad
	\HRule{\PropA,\neg\PropB}{\PropA\lor\PropB}{\HRuleName{\BKName}{18}}\quad
	\HRule{\PropA,\PropB}{\PropA\lor\PropB}{\HRuleName{\BKName}{19}}\\
	%\HRule{\PropA\to\PropB, \PropA}{\PropB}{\HRuleName{\BKName}{23}}\quad
	%\HRule{\PropA\lor\PropB, \PropC}{(\PropA\lor\neg\PropA)\to\PropC}{\HRuleName{\BKName}{24}}\quad
	%\HRule{\PropA\lor\PropB, \PropC}{\PropA\to\PropC}{\HRuleName{\BKName}{25}}\\
	%\HRule{\PropA\lor\PropB, (\PropA \lor \neg\PropA)\to\PropC}{\PropC}{\HRuleName{\BKName}{27}}\quad
 \HRule{\PropA\lor\PropB,\neg\PropA}{\PropB}{\HRuleName{\BKName}{20}}\quad
	\HRule{\PropA\lor(\PropB\lor\PropC)}{(\PropA\lor\PropB)\lor\PropC}{\HRuleName{\BKName}{21}}\quad
	\HRule{\PropA\lor\PropA}{\PropA}{\HRuleName{\BKName}{22}}\quad
	\HRule{\PropA\lor\PropB}{\PropB\lor\PropA}{\HRuleName{\BKName}{23}}\quad
	\HRule{\PropA\lor\PropB,\PropC}{\neg\PropA\lor\PropC}{\HRuleName{\BKName}{24}}\\
	%\HRule{\PropA\to(\PropB\to\PropC)}{(\PropA\land\PropB)\to\PropC}{\HRuleName{\BKName}{102}}\quad
	%\HRule{(\PropA\land\PropB)\to\PropC}{\PropA\to(\PropB\to\PropC)}{\HRuleName{\BKName}{103}}\\
	%\HRule{\PropA\to\PropB, \PropA\to\PropC}{\PropA\to(\PropB\land\PropC)}{\HRuleName{\BKName}{104}}\quad
\end{gather*}
plus the $\lor$-lifted versions of the above rules 
that satisfy the containment
condition (see Definition~\ref{def:containment_condition}) --- that is, all but $\HRuleName{\BKName}{1\star}$.
\end{definition}

\begin{remark}
%Adding the $\lor$-lifted versions of the  rules %displayed above substantially increases the 
%size of the proposed system. 
In classical logic, 
a rule can be derived from its $\lor$-lifted version
due to the presence of the rules 
$\HRule{\PropA\lor\PropA}{\PropA}{}$
and
$\HRule{\PropA}{\PropA\lor\PropB}{}$.
Since the latter is not sound in $\BKName$,
we needed to keep each $\HRuleName{\BKName}{i}$
and its $\lor$-lifted version
in the above calculus.
\end{remark}

Our first goal is to verify that the \SetFmla{} system just defined
is sound for $\BKName$.
This can be proved by showing
that each rule $\HRule{\FmSetC}{\FmD}{}$ of the system is \emph{sound for} $\SetFmlaRel_{\BKMat}$, i.e.~that $\FmSetC \SetFmlaRel_{\BKMat} \FmD$.
In this direction, we take advantage of the close
relationship between $\BKMat$ and classical logic described in
Theorem~\ref{the:bkcltheorem}.

\begin{lemma}
	\label{lem:soundness}
	$\BKVDash \; \subseteq \; \SetFmlaRel_{\BKMat}$.
\end{lemma}
\begin{proof}
	Note that $\HRuleName{\BKName}{1\star}$ is the only
	rule that does not satisfy the containment condition.
	Since it is impossible for an $\BKMat$-valuation 
	to satisfy both $\PropA$ and $\neg\PropA$, this
	rule is sound with respect to $\BKMat$.
	Because the other rules satisfy the containment condition
	and are all sound in classical logic, by
	Theorem~\ref{the:bkcltheorem}
	we have that they are also sound with respect to $\BKMat$.
\end{proof}

In what follows, we will abbreviate some 
\SetFmla{} derivations by composing rules of inference:
we write $\RuleA_1,\RuleA_2,\ldots,\RuleA_n$
to mean that we apply first rule $\RuleA_1$,
then $\RuleA_2$ considering the formula derived in the previous step,
then $\RuleA_3$ and so on.

\begin{proposition}
	\label{prop:basic_deriv_rules_bk}
	The following rules are derivable in $\ASymCalcName{\BKName}$:
	\begin{gather*}
	%	\HRule{\neg(\PropA\lor\PropB),\PropA}{\PropC}{\HRuleName{\BKName}{12\star}}\quad
	%	\HRule{\neg(\PropA\lor\PropB),\PropB}{\PropC}{\HRuleName{\BKName}{13\star}}\quad
		\HRule{\PropA\lor\PropB}{\PropA\to(\PropA\lor\PropB)}{\HRuleName{\BKName}{25}}\quad
		\HRule{\PropA\to\PropB, \PropA}{\PropB}{\HRuleName{\BKName}{26}}\quad
		%\HRule{\PropA\lor\PropB,\FmC}{\PropA\lor\FmC}{\HRuleName{\BKName}{32}}\quad
		\HRule{(\PropA\lor\PropB)\lor\PropC}{\PropA\lor(\PropB\lor\PropC)}{\HRuleName{\BKName}{27}}\quad
		%\HRule{\neg(\PropA\lor\PropB)}{\neg\PropA\land\neg\PropB}{\HRuleName{\BKName}{107}}\quad
		%\HRule{\neg\PropA\land\neg\PropB}{\neg(\PropA\lor\PropB)}{\HRuleName{\BKName}{108}}\\
		%\HRule{\PropA\to\PropB}{\neg(\PropB\land\neg\PropB)}{\HRuleName{\BKName}{30}}
		\HRule{\PropA\to\PropC,\PropB\to\PropC}{(\PropA\lor\PropB)\to\PropC}{\HRuleName{\BKName}{28}}\quad
		%\HRule{\PropA\to\PropB, \PropA\to\PropC}{\PropA\to(\PropB\land\PropC)}{\HRuleName{\BKName}{104}} (lifted version!)
	\end{gather*}
\end{proposition}
\begin{proof}{}
	Below we present the derivations of the above rules:
	\begin{itemize}
		%\item
		%$\HRule{\neg(\PropA\lor\PropB),\PropA}{\PropC}{\HRuleName{\BKName}{12\star}}$:
		%\begin{align*}
		%	n_1. \quad& \neg(\PropA\lor\PropB) & \text{Assumption}\\
		%	n_2. \quad& \PropA & \text{Assumption}\\
		%	k \SymbDef \max\{n_1,n_2\}+1. \quad& \neg\PropA & n_1, \HRuleName{\BKName}{15}\\
		%	k+1. \quad& \PropC & k, n_2,\HRuleName{\BKName}{1\star}\\
		%\end{align*}
%
		%\item
		%$\HRule{\neg(\PropA\lor\PropB),\PropB}{\PropC}{\HRuleName{\BKName}{13\star}}$: similar to the derivation of
		%$\HRuleName{\BKName}{12\star}$.
		\item $\HRule{\PropA\lor\PropB}{\PropA\to(\PropA\lor\PropB)}{\HRuleName{\BKName}{25}}$:
		\begin{align*}
			1. \quad& \PropA\lor\PropB & \text{Assumption}\\
			2. \quad& \PropA\lor\neg\PropA & 1, \HRuleName{\BKName}{15\star}\\
			3. \quad& (\PropA\lor\PropB)\lor(\PropA\lor\neg\PropA) & 2, \HRuleName{\BKName}{19}\\
			4. \quad& \neg\PropA\lor((\PropA\lor\PropB)\lor\PropA) & 3, \HRuleName{\BKName}{21},\HRuleName{\BKName}{23}\\
			5. \quad& \neg\PropA\lor((\PropA\lor\PropA)\lor\PropB) & 4, \HRuleName{\BKName}{23}^\lor,\HRuleName{\BKName}{21}^\lor\\
			6. \quad& \neg\PropA\lor(\PropA\lor\PropB) & 5, \HRuleName{\BKName}{23}^\lor, \HRuleName{\BKName}{22}^\lor, \HRuleName{\BKName}{23}^\lor\\
		\end{align*}
		\item $\HRule{\PropA\to\PropB, \PropA}{\PropB}{\HRuleName{\BKName}{26}}$: clearly from $\HRuleName{\BKName}{2}$ and $\HRuleName{\BKName}{20}$.
		\item $\HRule{(\PropA\lor\PropB)\lor\PropC}{\PropA\lor(\PropB\lor\PropC)}{\HRuleName{\BKName}{27}}$:
			clearly from $\HRuleName{\BKName}{21}$ and $\HRuleName{\BKName}{23}$.

		\item $\HRule{\PropA\to\PropC,\PropB\to\PropC}{(\PropA\lor\PropB)\to\PropC}{\HRuleName{\BKName}{28}}$:
		\begin{align*}
			1. \quad& \neg\PropA \lor \PropC & \text{Assumption}\\
			2. \quad& \neg\PropB \lor \PropC & \text{Assumption}\\
			3. \quad& \PropC\lor\neg\PropA & 1, \HRuleName{\BKName}{23}\\
			4. \quad& \PropC\lor\neg\PropB &2,\HRuleName{\BKName}{23}\\
			5. \quad& \PropC\lor\neg(\PropA\lor\PropB) &3, 4,\HRuleName{\BKName}{12}^\lor\\
			6. \quad& \neg(\PropA\lor\PropB)\lor\PropC &5,\HRuleName{\BKName}{23}\\
		\end{align*}
	\end{itemize}
\end{proof}

Now, should we also add as primitive rules the $\lor$-lifted versions of the
primitive $\lor$-lifted rules (and continue this {\it{ad infinitum}})? The following result 
shows that this is not necessary.

%{\color{red} Let us consider the abbreviation  $\varphi\to \psi=\neg %\varphi\lor \psi$.
%}

\begin{lemma}
	\label{lem:liftedderivableor}
	For every primitive rule $\RuleA$ of $\ASymCalcName{\BKName}$ but $\HRuleName{\BKName}{1\star}$,
	the $\lor$-lifted version of $\RuleA$ is derivable.
\end{lemma}
\begin{proof}
	Note that the $\lor$-lifted version of the
	rules depicted in Definition~\ref{def:bksetfmla},
	with the exception of $\HRuleName{\BKName}{1\star}$,
	are primitive in $\ASymCalcName{\BKName}$.
	Thus it remains to show that the $\lor$-lifted versions
	thereof are derivable in this system.
	Let
	$\HRule{\PropC \lor \FmA_1, \ldots, \PropC \lor \FmA_m}{\PropC\lor\FmB}{\RuleA^{\lor}}$ be
	the $\lor$-lifted version of
	$\HRule{\FmA_1, \ldots, \FmA_m}{\FmB}{\RuleA}$, this one being any of the
	primitive rules of $\ASymCalcName{\BKName}$ but $\HRuleName{\BKName}{1\star}$.
	Below we show that
	$\HRule{\PropD\lor(\PropC \lor \FmA_1), \ldots, \PropD\lor(\PropC \lor \FmA_m)}{\PropD\lor(\PropC\lor\FmB)}{\RuleA^{\lor\lor}}$
	is derivable in $\ASymCalcName{\BKName}$:
	\begin{align*}
		1.\quad & \PropD\lor(\PropC \lor \FmA_1) & \text{Assumption}\\
		    & \quad\vdots\\
		m.\quad & \PropD\lor(\PropC \lor \FmA_m) & \text{Assumption}\\
		m + 1.\quad & (\PropD\lor\PropC) \lor \FmA_1 & 1,\HRuleName{\BKName}{21}\\
		    & \quad\vdots\\
		2m.\quad & (\PropD\lor\PropC) \lor \FmA_m & m,\HRuleName{\BKName}{21}\\
		2m+1.\quad & (\PropD\lor\PropC) \lor \FmB & m+1,\ldots,2m,\RuleA^\lor\\
		2m+2.\quad & \PropD\lor(\PropC \lor \FmB) & 2m+1,\HRuleName{\BKName}{27}\\
	\end{align*}
\end{proof}

With the above, we also obtain the following result, which
will be useful to abbreviate some of the upcoming proofs.

\begin{corollary}
	\label{coro:liftderivableimp}
	For every primitive rule $\RuleA$ of $\ASymCalcName{\BKName}$ but $\HRuleName{\BKName}{1\star}$,
	the $\to$-lifted version of $\RuleA$ is derivable.
\end{corollary}
\begin{proof}
Let $\HRule{\FmA_1, \ldots, \FmA_m}{\FmB}{\RuleA}$ be a primitive rule of
$\ASymCalcName{\BKName}$ but $\HRuleName{\BKName}{1\star}$.
Then, from $\neg\PropC \lor \FmA_1, \ldots, \neg\PropC \lor \FmA_m$,
we derive, in view of Lemma~\ref{lem:liftedderivableor},
$\neg\PropC \lor \FmB$, and we are done.
\end{proof}

These two results extend easily to
rules that can be proved derivable in $\ASymCalcName{\BKName}$
without the use $\HRuleName{\BKName}{1\star}$.

\begin{corollary}
	\label{coro:lift_derivable_rules}
	Let $\RuleA$ be a derivable rule of $\ASymCalcName{\BKName}$
	having a proof that does not use $\HRuleName{\BKName}{1\star}$.
	Then $\RuleA^\lor$ and $\RuleA^\to$ are derivable
	as well.
\end{corollary}
\begin{proof}
	By induction on the length of the proof of $\RuleA$
	in $\ASymCalcName{\BKName}$ (one that does not employ $\HRuleName{\BKName}{1\star}$), applying essentially Lemma~\ref{lem:liftedderivableor}
	and Corollary~\ref{coro:liftderivableimp}.
\end{proof}

Even though $\BKName$ does not admit a deduction theorem
in the usual sense (see Theorem~\ref{fact:no-disj-no-ded-theo}), the following result provides
analogous deduction theorems that will be enough for our purposes.

\begin{proposition}
	\label{prop:dedtheclassphi}
    Let $\FmD \in \Set{\FmA,\FmB} \subseteq \LangSetInfec$
    and let $\TreeA$ be a proof in $\ASymCalcName{\BKName}$ witnessing
	that $\FmSetA, \FmA\lor\FmB, \FmD \SetFmlaRel_{\ASymCalcName{\BKName}} \FmC$.
 %for $\FmD$ being either $\FmA$ or $\FmB$.
	If the rule $\HRuleName{\BKName\star}{1}$ was not applied
	in $\TreeA$, then $\FmSetA,\FmA\lor\FmB \SetFmlaRel_{\ASymCalcName{\BKName}} \FmD \to \FmC$.
\end{proposition}
\begin{proof}
    
    Let us first consider the case $\FmD = \FmA$.
    Suppose that $\TreeA$ is $\FmC_1,\ldots,\FmC_n = \FmC$.
	We will prove that $P(j) \SymbDef \FmSetA,\FmA\lor\FmB \SetFmlaRel_{\ASymCalcName{\BKName}} \FmA \to \FmC_j$
	for all $1 \leq j \leq n$, using strong induction on $j$.
	For the base case $j=1$, we have that $\FmC_1 \in \FmSetA \cup \{\FmA\lor\FmB,\FmA\}$, leading to the following cases:
	\begin{enumerate}
		\item if $\FmC_1 \in \FmSetA$, use $\HRuleName{\BKName}{24}$.
		\item if $\FmC_1 = \FmA\lor\FmB$, use $\HRuleName{\BKName}{25}$.
		\item if $\FmC_1 = \FmA$, use $\HRuleName{\BKName}{15\star}$.
	\end{enumerate}
	Suppose now that (IH): $P(j)$ holds for all $j < k$. We want to prove $P(k)$.
	The cases when $\FmC_k \in \FmSetA \cup \{\FmA\lor\FmB,\FmA\}$ are as in the base case.
	We have to consider then $\FmC_k$ resulting from applications of the rules of
	the system, except for $\HRuleName{\BKName}{1\star}$.
	Assume that $\FmC_k$ resulted from an application of
	an $m$-ary rule $\RuleA$ using formulas $\FmC_{k_1},\ldots,\FmC_{k_m}$ as premises,
	which must have appeared previously in the proof.
	By (IH), then, we have
	$\FmSetA,\FmA\lor\FmB \SetFmlaRel_{\ASymCalcName{\BKName}} \FmA \to \FmC_{k_i}$
	for each $1 \leq i \leq m$.
	By Corollary~\ref{coro:liftderivableimp}, then, we have
	$\FmSetA,\FmA\lor\FmB \SetFmlaRel_{\ASymCalcName{\BKName}} \FmA \to \FmC_{k}$.
	In particular, for $k=n$, we obtain 
	$\FmSetA,\FmA\lor\FmB \SetFmlaRel_{\ASymCalcName{\BKName}} \FmA \to \FmC$, as desired.
 The case $\FmD = \FmB$ 
 follows easily by commutativity of $\lor$ and the 
 case $\FmD = \FmA$ just proved.
\end{proof}

%\begin{proposition}
%	\label{prop:dedtheclasspsi}
%	Let $\FmC_1,\ldots,\FmC_n$ be a proof in $\ASymCalcName{\BKName}$ witnessing
%	that $\FmSetA, \FmA\lor\FmB, \FmB \SetFmlaRel_{\ASymCalcName{\BKName}} \FmC$.
%	If the rule $\HRuleName{\BKName}{1}$ was not applied
%	in such proof, then $\FmSetA,\FmA\lor\FmB \SetFmlaRel_{\ASymCalcName{\BKName}} \FmB \to \FmC$.
%\end{proposition}
%\begin{proof}
%\end{proof}

With this deduction theorem, we can derive some rules
more easily, as the next result shows.
\begin{proposition}
	The following rules are derivable in $\ASymCalcName{\BKName}$:
	\begin{gather*}
		\HRule{\neg\PropA\lor\neg\PropB}{\neg(\PropA\land\PropB)}{\HRuleName{\BKName}{29}}\quad
		\HRule{\PropA\to\PropB}{\neg(\PropB\land\neg\PropB)}{\HRuleName{\BKName}{30}}
	\end{gather*}
\end{proposition}
\begin{proof}
    We present below the derivations.
    \mbox{}\\*	
    \begin{itemize}
		\item
		$\HRule{\neg\PropA\lor\neg\PropB}{\neg(\PropA\land\PropB)}{\HRuleName{\BKName}{29}}$:
		first of all, we prove that $\neg\PropA\lor\neg\PropB,\neg\PropB\lor\neg\neg\PropB, \neg\PropA \BKVDash \neg(\PropA\land\PropB)$:
		\begin{align*}
			1. \quad& \neg\PropA & \text{Assumption}\\
			2. \quad& \neg\PropB\lor\neg\neg\PropB & \text{Assumption}\\
			3. \quad& \neg\neg\PropB \lor \neg\PropA& 1, 2, \HRuleName{\BKName}{24}\\
			4. \quad& \neg\neg\PropB \lor \neg\PropB & 2, \HRuleName{\BKName}{23}&\\
			5. \quad& \neg\neg\PropB \lor \neg(\PropA\land\PropB)& 3,4,\HRuleName{\BKName}{5}^\lor\\
			6. \quad& \neg\neg\PropB \lor \neg\neg\neg\PropB & 2, \HRuleName{\BKName}{16\star}\\
			7. \quad& \neg\neg\neg\PropB \lor \neg\neg\PropB & 6, \HRuleName{\BKName}{23}\\
			8. \quad& \neg\neg\neg\PropB \lor \PropB & 7, \HRuleName{\BKName}{3}^\lor\\
			9. \quad& \neg\neg\neg\PropB \lor \neg\PropA & 1,3,\HRuleName{\BKName}{24}\\
			10.\quad& \neg\neg\neg\PropB \lor \neg(\PropA\land\PropB) & 8,9,\HRuleName{\BKName}{6}^\lor\\
			11.\quad& \neg\PropB \to \neg(\PropA\land\PropB) & 5, \text{Def. of $\to$}\\
			12.\quad& \neg\neg\PropB \to \neg(\PropA\land\PropB) & 10, \text{Def. of $\to$}\\
			13.\quad& (\neg\PropB\lor\neg\neg\PropB) \to \neg(\PropA\land\PropB) & 11, 12, \HRuleName{\BKName}{28}\\
			14.\quad& \neg(\PropA\land\PropB) & 2,13,\HRuleName{\BKName}{26}\\
		\end{align*}
		\noindent Similarly, we can show that 
		$\neg\PropA\lor\neg\PropB,\neg\PropA\lor\neg\neg\PropA, \neg\PropB \BKVDash \neg(\PropA\land\PropB)$,
		also without using $\HRuleName{\BKName}{1\star}$.
		Since $\HRuleName{\BKName}{1\star}$ was not employed in such derivations,
		we have
		$\neg\PropA\lor\neg\PropB,\neg\PropB\lor\neg\neg\PropB \BKVDash \neg\PropA \to \neg(\PropA\land\PropB)$
		and 
		$\neg\PropA\lor\neg\PropB,\neg\PropA\lor\neg\neg\PropA \BKVDash \neg\PropB \to \neg(\PropA\land\PropB)$,
		by Proposition~\ref{prop:dedtheclassphi}.
		Since 
		$\neg\PropA\lor\neg\PropB \BKVDash \neg\PropA\lor\neg\neg\PropA$ (by $\HRuleName{\BKName}{15\star}$)
		and
		$\neg\PropA\lor\neg\PropB \BKVDash \neg\PropB\lor\neg\neg\PropB$ (by $\HRuleName{\BKName}{16\star}$),
		by \REV{2.9}{transitivity of $\BKVDash$ (the \SetFmla{} version of \ref{prop:CRSSPropC})}, we have
		$\neg\PropA\lor\neg\PropB \BKVDash \neg\PropA \to \neg(\PropA\land\PropB)$
		and 
		$\neg\PropA\lor\neg\PropB \BKVDash \neg\PropB \to \neg(\PropA\land\PropB)$.
		By $\HRuleName{\BKName}{28}$, then,
		$\neg\PropA\lor\neg\PropB \BKVDash (\neg\PropA \lor \neg\PropB) \to \neg(\PropA\land\PropB)$.
		Finally, by $\HRuleName{\BKName}{26}$ (modus ponens), we obtain
		$\neg\PropA\lor\neg\PropB \BKVDash \neg(\PropA\land\PropB)$.
		\item $\HRule{\PropA\to\PropB}{\neg(\PropB\land\neg\PropB)}{\HRuleName{\BKName}{30}}$:
			clearly from $\HRuleName{\BKName}{16\star}$ and $\HRuleName{\BKName}{29}$.
	\end{itemize}
\end{proof}

In the negation fragment of classical logic (call it $\CLName_{\neg}$)
we have the following deduction theorem:
if $\FmSetA, \FmA \SetFmlaRel_{\CLName_{\neg}} \neg\FmA$, then
	$\FmSetA \SetFmlaRel_{\CLName_{\neg}} \neg\FmA$.
 Similarly to what we did in Proposition~\ref{prop:dedtheclassphi}, we show
 now that this result also holds for $\BKName$
 provided $\FmA \lor \FmB$ is present in the context $\FmSetA$.
 In this case, however, we do not need
 to impose any restriction on the
 rules applied in the derivations witnessing
 the consecution in the assumption.
 This result will also be useful for proving the desired completeness result
 for $\ASymCalcName{\BKName}$.

\begin{proposition}
	\label{prop:negdedthe}
    Let $\FmD \in \Set{\FmA,\FmB} \subseteq \LangSetInfec$.
    If $\FmSetA, \FmA\lor\FmB, \FmD \SetFmlaRel_{\ASymCalcName{\BKName}} \neg\FmD$, then
	$\FmSetA, \FmA\lor\FmB \SetFmlaRel_{\ASymCalcName{\BKName}} \neg\FmD$.
	%If $\FmSetA, \FmA\lor\FmB, \FmA \SetFmlaRel_{\ASymCalcName{\BKName}} \neg\FmA$, then
	%$\FmSetA, \FmA\lor\FmB \SetFmlaRel_{\ASymCalcName{\BKName}} \neg\FmA$.
\end{proposition}
\begin{proof}
We begin with the case $\FmD = \FmA$.
Let $\TreeA = \FmC_1,\ldots,\FmC_n$ be a proof witnessing that
$\FmSetA, \FmA\lor\FmB, \FmA \SetFmlaRel_{\ASymCalcName{\BKName}} \neg\FmA$.
In case no application of $\HRuleName{\BKName}{1\star}$ is used in $\TreeA$,
we have 
$\FmSetA, \FmA\lor\FmB \SetFmlaRel_{\ASymCalcName{\BKName}} \FmA \to \neg\FmA \SetFmlaRel_{\ASymCalcName{\BKName}}  \neg\FmA \lor \neg\FmA \SetFmlaRel_{\ASymCalcName{\BKName}} \neg\FmA$ by Proposition~\ref{prop:dedtheclassphi} and $\HRuleName{\BKName}{22}$, as desired.
On the other hand, suppose that $\FmC_k$ was the formula produced by the first application
	of $\HRuleName{\BKName}{1\star}$. Then $k > 2$ and there are
	$\FmC_{m_1}$ and $\FmC_{m_2} = \neg\FmC_{m_1}$, with $m_1,m_2 < k$, such that
	$\FmSetA,\FmA\lor\FmB, \FmA \SetFmlaRel_{\ASymCalcName{\BKName}} \FmC_{m_1}$ and 
	$\FmSetA,\FmA\lor\FmB, \FmA \SetFmlaRel_{\ASymCalcName{\BKName}} \neg\FmC_{m_1}$.
	By Proposition~\ref{prop:dedtheclassphi}, then, we have
	(a): $\FmSetA,\FmA\lor\FmB \SetFmlaRel_{\ASymCalcName{\BKName}} \FmA \to \FmC_{m_1}$,
	and
	$\FmSetA,\FmA\lor\FmB \SetFmlaRel_{\ASymCalcName{\BKName}} \FmA \to \neg\FmC_{m_1}$,
	and so, by Corollary~\ref{coro:liftderivableimp},
	$\FmSetA,\FmA\lor\FmB \SetFmlaRel_{\ASymCalcName{\BKName}} \FmA \to (\FmC_{m_1} \land \neg\FmC_{m_1})$,
	that is,
	$\FmSetA,\FmA\lor\FmB \SetFmlaRel_{\ASymCalcName{\BKName}} \neg\FmA \lor (\FmC_{m_1} \land \neg\FmC_{m_1}) \SetFmlaRel_{\ASymCalcName{\BKName}}
	 (\FmC_{m_1} \land \neg\FmC_{m_1})\lor\neg\FmA$.
	 But $\FmSetA,\FmA\lor\FmB \SetFmlaRel_{\ASymCalcName{\BKName}} \neg(\FmC_{m_1} \land \neg\FmC_{m_1})$
	 by (a) and $\HRuleName{\BKName}{30}$,
	 and thus 
	 $\FmSetA,\FmA\lor\FmB \SetFmlaRel_{\ASymCalcName{\BKName}} \neg\FmA$
  by $\HRuleName{\BKName}{20}$.
  Now, for $\FmD = \FmB$,
  we have that from
$\FmSetA, \FmA\lor\FmB, \FmB \SetFmlaRel_{\ASymCalcName{\BKName}} \neg\FmB$
and the rule $\HRuleName{\BKName}{23}$, we get
$\FmSetA, \FmB\lor\FmA, \FmB \SetFmlaRel_{\ASymCalcName{\BKName}} \neg\FmB$.
By Proposition~\ref{prop:negdedthe}, we get
$\FmSetA, \FmB\lor\FmA \SetFmlaRel_{\ASymCalcName{\BKName}} \neg\FmB$
and, again by  $\HRuleName{\BKName}{23}$, we have
$\FmSetA, \FmA\lor\FmB \SetFmlaRel_{\ASymCalcName{\BKName}} \neg\FmB$.
\end{proof}

%\begin{proposition}
%	\label{prop:negdedthepsi}
%	If $\FmSetA, \FmA\lor\FmB, \FmB \SetFmlaRel_{\ASymCalcName{\BKName}} \neg\FmB$, then
%	$\FmSetA, \FmA\lor\FmB \SetFmlaRel_{\ASymCalcName{\BKName}} \neg\FmB$.
%\end{proposition}
%\begin{proof}
%	From
%$\FmSetA, \FmA\lor\FmB, \FmB \SetFmlaRel_{\ASymCalcName{\BKName}} \neg\FmB$
%and the rule $\HRuleName{\BKName}{23}$, we get
%$\FmSetA, \FmB\lor\FmA, \FmB \SetFmlaRel_{\ASymCalcName{\BKName}} \neg\FmB$.
%By Proposition~\ref{prop:negdedthe}, we get
%$\FmSetA, \FmB\lor\FmA \SetFmlaRel_{\ASymCalcName{\BKName}} \neg\FmB$
%and, again by  $\HRuleName{\BKName}{23}$, we have
%$\FmSetA, \FmA\lor\FmB \SetFmlaRel_{\ASymCalcName{\BKName}} \neg\FmB$.
%\end{proof}

A consequence of the previous result
is the following.

\begin{proposition}
	\label{prop:dedthewithexppsi}
 Let $\FmD_1,\FmD_2 \in \Set{\FmA,\FmB} \subseteq \LangSetInfec$
 with $\FmD_1 \neq \FmD_2$,
 and $\FmC_1,\ldots,\FmC_n$ be a proof in $\ASymCalcName{\BKName}$ witnessing
	that $\FmSetA, \FmA\lor\FmB, \FmD_1 \SetFmlaRel_{\ASymCalcName{\BKName}} \FmC$.
	If the rule $\HRuleName{\BKName}{1\star}$ was applied
	in such proof, then $\FmSetA, \FmA\lor\FmB \SetFmlaRel_{\ASymCalcName{\BKName}} \FmD_2$.
	%Let $\FmC_1,\ldots,\FmC_n$ be a proof in $\ASymCalcName{\BKName}$ witnessing
	%that $\FmSetA, \FmA\lor\FmB, \FmA \SetFmlaRel_{\ASymCalcName{\BKName}} \FmC$.
	%If the rule $\HRuleName{\BKName}{1\star}$ was applied
	%in such proof, then $\FmSetA, \FmA\lor\FmB \SetFmlaRel_{\ASymCalcName{\BKName}} \FmB$.
\end{proposition}
\begin{proof}
    We will prove the case $\FmD_1 = \FmA$
    and the other will be analogous.
	Suppose that $\FmC_k$ was the formula produced by the first application
	of $\HRuleName{\BKName}{1\star}$. Then $k > 2$ and there are
	$\FmC_{m_1}$ and $\FmC_{m_2} = \neg\FmC_{m_1}$, with $m_1,m_2 < k$, such that
	$\FmSetA,\FmA\lor\FmB, \FmA \SetFmlaRel_{\ASymCalcName{\BKName}} \FmC_{m_1}$ and 
	$\FmSetA,\FmA\lor\FmB, \FmA \SetFmlaRel_{\ASymCalcName{\BKName}} \neg\FmC_{m_1}$.
	But then, by $\HRuleName{\BKName}{1\star}$,
	we have
	$\FmSetA,\FmA\lor\FmB, \FmA \SetFmlaRel_{\ASymCalcName{\BKName}} \neg\FmA$.
	By Proposition~\ref{prop:negdedthe}, then, we have
	$\FmSetA,\FmA\lor\FmB \SetFmlaRel_{\ASymCalcName{\BKName}} \neg\FmA$,
	and then
	$\FmSetA,\FmA\lor\FmB \SetFmlaRel_{\ASymCalcName{\BKName}} \FmB$
	by $\HRuleName{\BKName}{20}$.
	%Since all formulas prior to $\FmC_k$ were derived by rules other than
	%$\HRuleName{\BKName}{1}$,
	%by Proposition~\ref{prop:dedtheclass},
	%we have
	%$\FmSetA,\FmA\lor\FmB \vdash \FmA\to\FmC_{m_1}$ and 
	%$\FmSetA,\FmA\lor\FmB \vdash \FmA\to\neg\FmC_{m_1}$,
	%and then
	%$\FmSetA,\FmA\lor\FmB \vdash \FmA\to(\FmC_{m_1} \land \neg\FmC_{m_1})$.
\end{proof}

%\begin{proposition}
%	\label{prop:dedthewithexpphi}
%	Let $\FmC_1,\ldots,\FmC_n$ be a proof in $\ASymCalcName{\BKName}$ witnessing
%	that $\FmSetA, \FmA\lor\FmB, \FmB \SetFmlaRel_{\ASymCalcName{\BKName}} \FmC$.
%	If the rule $\HRuleName{\BKName}{1\star}$ was applied
%	in such proof, then $\FmSetA, \FmA\lor\FmB \SetFmlaRel_{\ASymCalcName{\BKName}} \FmA$.
%\end{proposition}
%\begin{proof}
%	Similar to the previous proof.
%\end{proof}

In Proposition~\ref{fact:no-disj-no-ded-theo}, we proved that
$\BKName$ does not allow to express
a connective satisfying (\ref{eq:disj-prop}).
Nevertheless, we now show that having
$\FmA\lor\FmB$ in the context is also enough
to recover this result.

\begin{lemma}
	\label{lem:disjunction}
	For all $\FmSetA, \Set{\FmA,\FmB,\FmC} \subseteq \LangSetInfec$, we have
			%if
			$\FmSetA,\FmA\lor\FmB,\FmA \vdash_{\ASymCalcName{\BKName}}\FmC$
			and
			$\FmSetA,\FmA\lor\FmB,\FmB \vdash_{\ASymCalcName{\BKName}}\FmC$,
			if, and only if,
			$\FmSetA,\FmA\lor\FmB \vdash_{\ASymCalcName{\BKName}}\FmC$.
\end{lemma}
\begin{proof}
The right-to-left direction is obvious by
reflexivity of the consequence relation.
For the left-to-right direction,
suppose that $\TreeA_1$ and $\TreeA_2$ are witnesses of
$\FmSetA,\FmA\lor\FmB,\FmA \vdash_{\ASymCalcName{\BKName}}\FmC$
and
$\FmSetA,\FmA\lor\FmB,\FmB \vdash_{\ASymCalcName{\BKName}}\FmC$,
respectively.
Consider the following cases:
\begin{enumerate}
	\item In both there are no applications of $\HRuleName{\BKName}{1\star}$:
		then, by Proposition~\ref{prop:dedtheclassphi},
		we have
		$\FmSetA,\FmA\lor\FmB \vdash_{\ASymCalcName{\BKName}}\FmA\to\FmC$
		and
		$\FmSetA,\FmA\lor\FmB\vdash_{\ASymCalcName{\BKName}}\FmB\to\FmC$.
		Thus
		$\FmSetA,\FmA\lor\FmB\vdash_{\ASymCalcName{\BKName}}(\FmA\lor\FmB)\to\FmC$, by $\HRuleName{\BKName}{28}$,
		and
		$\FmSetA,\FmA\lor\FmB\vdash_{\ASymCalcName{\BKName}}\FmC$, by $\HRuleName{\BKName}{26}$.
	\item If there is an application of $\HRuleName{\BKName}{1\star}$ in $\TreeA_1$:
		then, by Proposition~\ref{prop:dedthewithexppsi},
		we have
		$\FmSetA,\FmA\lor\FmB\vdash_{\ASymCalcName{\BKName}}\FmB$.
		Then, by transitivity
		considering 
			$\FmSetA,\FmA\lor\FmB,\FmB \vdash_{\ASymCalcName{\BKName}}\FmC$,
			we obtain
			$\FmSetA,\FmA\lor\FmB\vdash_{\ASymCalcName{\BKName}}\FmC$.
	\item If there is an application of $\HRuleName{\BKName}{1\star}$ in $\TreeA_2$: similar to the previous case.
\end{enumerate}
\end{proof}

Finally we get to the desired axiomatization result.

\begin{theorem}
    $\SetFmlaRel_{\ASymCalcName{\BKName}} \; = \; \SetFmlaRel_{\BKName}$.
	%$\FmSetA \vdash_{\ASymCalcName{\BKName}}\FmB$	
	%if, and only if,
	%$\FmSetA \SetSetRel_{\SymCalcName{\BKName\star}}{\Set{\FmB}}$.
\end{theorem}
\begin{proof}
We will show that
$\FmSetA \vdash_{\ASymCalcName{\BKName}}\FmB$	
	if, and only if,
	$\FmSetA \SetSetRel_{\SymCalcName{\BKName\star}}{\Set{\FmB}}$.
	The left-to-right direction easily follows,
	since every rule of $\ASymCalcName{\BKName}$
	is sound with respect to the matrix of $\BKName$ by Lemma~\ref{lem:soundness},
	and thus derivable in $\SymCalcName{\BKName\star}$.
	From the right to the left,
	we will show by induction on the structure of derivations in $\SymCalcName{\BKName\star}$
	that $P(\TreeA)$: if $\TreeA$ witnesses that $\FmSetA \SetSetRel_{\SymCalcName{\BKName\star}}{\Set{\FmB}}$, 
	then there is a proof in $\ASymCalcName{\BKName}$ bearing witness to $\FmSetA\BKVDash\FmB$.
	In the base case, $\TreeA$ has a single node, meaning that $\FmB \in \FmSetA$, and we are done by
	reflexivity of $\SetFmlaRel_{\ASymCalcName{\BKName}}$.
	In the inductive step, we assume $P(\TreeA')$ for each subtree $\TreeA'$ of $\TreeA$ 
	and consider $\TreeA$ resulting from an application of the rules of $\SymCalcName{\BKName\star}$.
	Let us consider three cases:
	\begin{enumerate}
		\item $\TreeA$ results from a rule that is derivable in $\ASymCalcName{\BKName}$: here, there is
			nothing to do, as the same rule may be applied to produce the desired derivation.
		\item $\TreeA$ results from an application of $\HRuleSSName{\BKName}{1}$: use $\HRuleName{\BKName}{1\star}$ instead.
		\item $\TreeA$ results from an application of $\HRuleSSName{\BKName}{20}$:
			if the root of $\TreeA$ is labelled with $\FmSetC$, then
			$\FmC\lor\FmD \in \FmSetC$, 
			and we have, by the induction hypothesis,
			(a): $\FmSetC,\FmC\lor\FmD,\FmC\SetFmlaRel_{\ASymCalcName{\BKName}}\FmB$ and
			(b): $\FmSetC,\FmC\lor\FmD,\FmD\SetFmlaRel_{\ASymCalcName{\BKName}}\FmB$.
			By Lemma~\ref{lem:disjunction}, then,
			we obtain the desired result.
			%By Proposition~\ref{prop:dedthe}, we have
			%(a1): $\FmSetC,\FmC\lor\FmD\SetFmlaRel_{\ASymCalcName{\BKName}}\FmC\to\FmB$ or
			%(a2): $\FmSetC,\FmC\lor\FmD\SetFmlaRel_{\ASymCalcName{\BKName}}(\FmC\lor\neg\FmC)\to\FmB$
			%and
			%(b1): $\FmSetC,\FmC\lor\FmD\SetFmlaRel_{\ASymCalcName{\BKName}}\FmD\to\FmB$ or
			%(b2): $\FmSetC,\FmC\lor\FmD\SetFmlaRel_{\ASymCalcName{\BKName}}(\FmD\lor\neg\FmD)\to\FmB$.
			%In case (a2) or (b2), use rules 
			%$\HRuleName{\BKName}{17\star}$ or
			%$\HRuleName{\BKName}{18\star}$, and
			%$\HRuleName{\BKName}{23}$.
			%In case (a1) and (b1), 
			%we have $\FmSetC,\FmC\lor\FmD\SetFmlaRel_{\ASymCalcName{\BKName}}(\FmC\lor\FmD)\to\FmB$,
			%by $\HRuleName{\BKName}{26}$. Then, by
			%$\HRuleName{\BKName}{23}$, we have
			%$\FmSetC,\FmC\lor\FmD\SetFmlaRel_{\ASymCalcName{\BKName}}\FmB$,
			%and, since $\FmC\lor\FmD \in \FmSetC$,
			%$\FmSetC\SetFmlaRel_{\ASymCalcName{\BKName}}\FmB$.
	\end{enumerate}
\end{proof}

\section{Concluding remarks}\label{sec:final_remarks}

Taking stock of what we  achieved in
 the previous sections, 
we highlight that  we have settled  fundamental questions
regarding $\BKName$ and $\PWKName$,  two logics that are among the main subjects of this Special Issue.
We also wish to mention an interesting corollary of our results, namely
that 
some finite subset of the axioms employed in the papers~\cite{bonzio2021plonka,Bonzio2021-BONCLA}
must already suffice to axiomatize each of the two logics. We leave this observation as a suggestion for 
future developments. 

Besides the intrinsic interest in the results established above,
the present paper may also be seen as another illustration of 
the differences in expressive power among the various 
available proof-theoretic formalisms in logic, and in particular between
$\SetSet{}$ over $\SetFmla{}$ H-systems. The latter are obviously
less expressive than the former --  even weaker if compared to sequent systems -- 
even though they afford more fine-grained tools for comparing and also for combining logics (in particular when one wishes to introduce
the least possible interactions),
as recent results amply demonstrate~\cite{Fibring1,Fibring2}.

Another direction for future research worth mentioning is the study of these and other logics associated to the algebra  $\AlgInfec$ (and other three-valued algebras)
in the setting of different kinds of H-systems. In particular,
 a two-dimensional version of $\SetSet{}$ H-systems~\cite{greati2021,greati2022}, whose induced logics are the so-called  \emph{{\normalfont\textsf{B}}-consequence relations}~\cite{blasio2017},
 may be employed as a uniform setting for investigating
\emph{pure consequence relations} (like $\BKName$ and $\PWKName$), their intersection (\emph{order-theoretic consequence relations}) and \emph{mixed consequence relations} (we use here the terminology of~\cite{chemla2017}),
the latter being non-Tarskian consequence relations
(lacking either reflexivity~\cite{malinowski1994} or transitivity~\cite{frankowski2004form}). 

Not only can 
%As if it were not enough that 
a
two-dimensional logic  express all of these very different notions of logics in the same logical environment: we also
have that it has a neat analytic two-dimensional axiomatization.
That is, this two-dimensional logic has not only great theoretical value due to its expressiveness, but also constitutes an important tool for using the above-mentioned logics and studying their properties.

\subsubsection*{Acknowledgments}
Vitor Greati acknowledges support from the FWF project P33548.
Sérgio Marcelino's research was done under the scope of project FCT/MCTES through national funds and when applicable co-funded by EU under the project UIDB/50008/2020.

%Bibliography
\bibliographystyle{abbrvnat}  
\bibliography{sn-bibliography}

\begin{thebibliography}{35}
\providecommand{\natexlab}[1]{#1}
\providecommand{\url}[1]{\texttt{#1}}
\expandafter\ifx\csname urlstyle\endcsname\relax
  \providecommand{\doi}[1]{doi: #1}\else
  \providecommand{\doi}{doi: \begingroup \urlstyle{rm}\Url}\fi

\bibitem[Avron and Zamansky(2011)]{avron2011}
A.~Avron and A.~Zamansky.
\newblock Non-deterministic semantics for logical systems.
\newblock In D.~M. Gabbay and F.~Guenthner, editors, \emph{Handbook of Philosophical Logic: Volume 16}, pages 227--304. Springer, Dordrecht, 2011.
\newblock \doi{10.1007/978-94-007-0479-4_4}.

\bibitem[Baaz et~al.(1996)Baaz, Ferm{\"u}ller, Salzer, and Zach]{baaz1996}
M.~Baaz, C.~G. Ferm{\"u}ller, G.~Salzer, and R.~Zach.
\newblock Multlog 1.0: Towards an expert system for many-valued logics.
\newblock In M.~A. McRobbie and J.~K. Slaney, editors, \emph{Automated Deduction --- CADE13}, pages 226--230, Berlin, Heidelberg, 1996. Springer Berlin Heidelberg.
\newblock ISBN 978-3-540-68687-3.
\newblock \doi{10.1007/3-540-61511-3_84}.

\bibitem[Baaz et~al.(2013)Baaz, Lahav, and Zamansky]{baaz2013}
M.~Baaz, O.~Lahav, and A.~Zamansky.
\newblock Finite-valued semantics for canonical labelled calculi.
\newblock \emph{Journal of Automated Reasoning}, 51, 2013.
\newblock \doi{10.1007/s10817-013-9273-x}.

\bibitem[Belikov(2021)]{belikov2021}
A.~Belikov.
\newblock {On bivalent semantics and natural deduction for some infectious logics}.
\newblock \emph{Logic Journal of the IGPL}, 30\penalty0 (1):\penalty0 186--210, 02 2021.
\newblock ISSN 1367-0751.
\newblock \doi{10.1093/jigpal/jzaa071}.

\bibitem[Bochvar and Bergmann(1981)]{bochvar1981}
D.~Bochvar and M.~Bergmann.
\newblock On a three-valued logical calculus and its application to the analysis of the paradoxes of the classical extended functional calculus.
\newblock \emph{History and Philosophy of Logic}, 2\penalty0 (1-2):\penalty0 87--112, 1981.
\newblock \doi{10.1080/01445348108837023}.

\bibitem[Bochvar(1938)]{bochvar1938original}
D.~A. Bochvar.
\newblock On a three-valued logical calculus and its application to the analysis of contradictions.
\newblock \emph{Rec. Math. Moscou, n. Ser.}, 4:\penalty0 287--308, 1938.

\bibitem[Bonzio and Baldi(2021)]{Bonzio2021-BONCLA}
S.~Bonzio and M.~P. Baldi.
\newblock Containment logics: algebraic completeness and axiomatization.
\newblock \emph{Studia Logica}, 109\penalty0 (5):\penalty0 969--994, 2021.
\newblock \doi{10.1007/s11225-020-09930-1}.

\bibitem[Bonzio et~al.(2017)Bonzio, Gil-F{\'{e}}rez, Paoli, and Peruzzi]{BonGilPaoPer17}
S.~Bonzio, J.~Gil-F{\'{e}}rez, F.~Paoli, and L.~Peruzzi.
\newblock {On {P}araconsistent {W}eak {K}leene logic: axiomatisation and algebraic analysis}.
\newblock \emph{Studia Logica}, 105\penalty0 (2):\penalty0 253--297, 2017.
\newblock ISSN 1572-8730.
\newblock \doi{10.1007/s11225-016-9689-5}.
\newblock URL \url{https://doi.org/10.1007/s11225-016-9689-5}.

\bibitem[Bonzio et~al.(2021)Bonzio, Moraschini, and Pra~Baldi]{bonzio2021plonka}
S.~Bonzio, T.~Moraschini, and M.~Pra~Baldi.
\newblock Logics of left variable inclusion and {P}{\l}onka sums of matrices.
\newblock \emph{Archive for Mathematical Logic}, 60\penalty0 (1):\penalty0 49--76, Feb 2021.
\newblock ISSN 1432-0665.
\newblock \doi{10.1007/s00153-020-00727-6}.
\newblock URL \url{https://doi.org/10.1007/s00153-020-00727-6}.

\bibitem[Bonzio et~al.(2022{\natexlab{a}})Bonzio, Paoli, and Baldi]{Bonzio2022book}
S.~Bonzio, F.~Paoli, and M.~P. Baldi.
\newblock \emph{Logics of Variable Inclusion}.
\newblock Springer International Publishing, Cham, 2022{\natexlab{a}}.
\newblock \doi{10.1007/978-3-031-04297-3}.

\bibitem[Bonzio et~al.(2022{\natexlab{b}})Bonzio, Paoli, and Pra~Baldi]{Bonzio2022}
S.~Bonzio, F.~Paoli, and M.~Pra~Baldi.
\newblock \emph{Paraconsistent Weak Kleene Logic}, pages 159--198.
\newblock Springer International Publishing, Cham, 2022{\natexlab{b}}.
\newblock ISBN 978-3-031-04297-3.
\newblock \doi{10.1007/978-3-031-04297-3_7}.

\bibitem[C.~Blasio(2017)]{blasio2017}
H.~W. C.~Blasio, J.~Marcos.
\newblock An inferentially many-valued two-dimensional notion of entailment.
\newblock \emph{Bulletin of the Section of Logic}, 46\penalty0 (3/4), 2017.
\newblock \doi{10.18778/0138-0680.46.3.4.05}.

\bibitem[Caleiro and Marcelino(2019)]{marcelino19woll}
C.~Caleiro and S.~Marcelino.
\newblock Analytic calculi for monadic {P}{N}matrices.
\newblock In R.~Iemhoff, M.~Moortgat, and R.~Queiroz, editors, \emph{Logic, Language, Information and Computation (WoLLIC 2019)}, volume 11541 of \emph{LNCS}, pages 84--98. Springer, Cham, 2019.
\newblock \doi{10.1007/978-3-662-59533-6_6}.

\bibitem[Caleiro and Marcelino(2023)]{Fibring2}
C.~Caleiro and S.~Marcelino.
\newblock Modular semantics for combined many-valued logics.
\newblock \emph{Notre Dame Journal of Formal Logic}, pages 1--54, 2023.
\newblock \doi{https://doi.org/10.1017/jsl.2023.22}.

\bibitem[Caleiro et~al.(2020)Caleiro, Marcelino, and Filipe]{caleiro2020infec}
C.~Caleiro, S.~Marcelino, and P.~Filipe.
\newblock Infectious semantics and analytic calculi for even more inclusion logics.
\newblock In \emph{2020 IEEE 50th International Symposium on Multiple-Valued Logic (ISMVL)}, pages 224--229, 2020.
\newblock \doi{10.1109/ISMVL49045.2020.000-1}.

\bibitem[Chemla et~al.(2017)Chemla, Égré, and Spector]{chemla2017}
E.~Chemla, P.~Égré, and B.~Spector.
\newblock {Characterizing logical consequence in many-valued logic}.
\newblock \emph{Journal of Logic and Computation}, 27\penalty0 (7):\penalty0 2193--2226, 03 2017.
\newblock ISSN 0955-792X.
\newblock \doi{10.1093/logcom/exx001}.
\newblock URL \url{https://doi.org/10.1093/logcom/exx001}.

\bibitem[Da~Ré et~al.(2018)Da~Ré, Pailos, and Szmuc]{pailos2018}
B.~Da~Ré, F.~Pailos, and D.~Szmuc.
\newblock {Theories of truth based on four-valued infectious logics}.
\newblock \emph{Logic Journal of the IGPL}, 28\penalty0 (5):\penalty0 712--746, 11 2018.
\newblock ISSN 1367-0751.
\newblock \doi{10.1093/jigpal/jzy057}.

\bibitem[Ferguson(2014)]{ferguson2014}
T.~Ferguson.
\newblock A computational interpretation of conceptivism.
\newblock \emph{Journal of Applied Non-Classical Logics}, 24\penalty0 (4):\penalty0 333--367, 2014.
\newblock \doi{10.1080/11663081.2014.980116}.

\bibitem[Frankowski(2004)]{frankowski2004form}
S.~Frankowski.
\newblock Formalization of a plausible inference.
\newblock \emph{Bulletin of the Section of Logic}, 33:\penalty0 41--52, 01 2004.

\bibitem[Greati(2022)]{greatimsc2021}
V.~Greati.
\newblock Hilbert-style formalism for two-dimensional notions of consequence.
\newblock Master's thesis, Universidade Federal do Rio Grande do Norte, 2022.
\newblock URL \url{https://repositorio.ufrn.br/handle/123456789/46792}.

\bibitem[Greati and Marcos(2022)]{greati2022}
V.~Greati and J.~Marcos.
\newblock Finite two-dimensional proof systems for non-finitely axiomatizable logics.
\newblock In J.~Blanchette, L.~Kov{\'a}cs, and D.~Pattinson, editors, \emph{Automated Reasoning}, pages 640--658, Cham, 2022. Springer International Publishing.
\newblock ISBN 978-3-031-10769-6.
\newblock \doi{10.1007/978-3-031-10769-6_37}.

\bibitem[Greati et~al.()Greati, Greco, Marcelino, Palmigiano, and Rivieccio]{3valbook}
V.~Greati, G.~Greco, S.~Marcelino, A.~Palmigiano, and U.~Rivieccio.
\newblock chapter Generating proof systems for three-valued propositional logics (submitted).

\bibitem[Greati et~al.(2021)Greati, Marcelino, and Marcos]{greati2021}
V.~Greati, S.~Marcelino, and J.~Marcos.
\newblock Proof search on bilateralist judgments over non-deterministic semantics.
\newblock In A.~Das and S.~Negri, editors, \emph{Automated Reasoning with Analytic Tableaux and Related Methods}, pages 129--146, Cham, 2021. Springer International Publishing.
\newblock ISBN 978-3-030-86059-2.
\newblock \doi{10.1007/978-3-030-86059-2_8}.

\bibitem[Halld\'{e}n(1949)]{hallden1949}
S.~Halld\'{e}n.
\newblock \emph{The Logic of Nonsense}.
\newblock Upsala Universitets Arsskrift, Uppsala, Sweden, 1949.

\bibitem[Kleene(1952)]{Kl50}
S.~C. Kleene.
\newblock \emph{Introduction to Metamathematics}.
\newblock North Holland, Princeton, NJ, USA, 1952.

\bibitem[Malinowski(1994)]{malinowski1994}
G.~Malinowski.
\newblock \emph{Inferential many-valuedness}, pages 75--84.
\newblock Springer Netherlands, Dordrecht, 1994.
\newblock ISBN 978-94-015-8273-5.
\newblock \doi{10.1007/978-94-015-8273-5_6}.

\bibitem[Marcelino and Caleiro(2017)]{Fibring1}
S.~Marcelino and C.~Caleiro.
\newblock Disjoint fibring of non-deterministic matrices.
\newblock In R.~de~Queiroz and J.~Kennedy, editors, \emph{Logic, Language, Information and Computation (WoLLIC 2017)}, volume 10388 of \emph{LNCS}, pages 242--255. Springer, Berlin, Heidelberg, 2017.
\newblock \doi{10.1007/978-3-662-55386-2_17}.

\bibitem[Marcelino and Caleiro(2019)]{marcelino19syn}
S.~Marcelino and C.~Caleiro.
\newblock Axiomatizing non-deterministic many-valued generalized consequence relations.
\newblock \emph{Synthese}, 198, 2019.
\newblock \doi{10.1007/s11229-019-02142-8}.

\bibitem[Negri et~al.(2001)Negri, von Plato, and Ranta]{negri2001}
S.~Negri, J.~von Plato, and A.~Ranta.
\newblock \emph{Structural Proof Theory}.
\newblock Cambridge University Press, Cambridge, 2001.
\newblock \doi{10.1017/CBO9780511527340}.

\bibitem[Paoli and Pra~Baldi(2020)]{Paoli2020}
F.~Paoli and M.~Pra~Baldi.
\newblock Proof theory of {P}araconsistent {W}eak {K}leene logic.
\newblock \emph{Studia Logica}, 108\penalty0 (4):\penalty0 779--802, Aug 2020.
\newblock ISSN 1572-8730.
\newblock \doi{10.1007/s11225-019-09876-z}.
\newblock URL \url{https://doi.org/10.1007/s11225-019-09876-z}.

\bibitem[Petrukhin(2017)]{Petrukhin2017}
Y.~Petrukhin.
\newblock Natural deduction for three-valued regular logics.
\newblock \emph{Logic and Logical Philosophy}, 26\penalty0 (2):\penalty0 197--206, 2017.

\bibitem[Prior(1967)]{prior1967}
A.~N. Prior.
\newblock \emph{Past, Present, and Future}.
\newblock Clarendon Press, Oxford, England, 1967.

\bibitem[Shoesmith and Smiley(1978)]{ss1978}
D.~J. Shoesmith and T.~J. Smiley.
\newblock \emph{Multiple-Conclusion Logic}.
\newblock Cambridge University Press, Cambridge, 1978.
\newblock \doi{10.1017/CBO9780511565687}.

\bibitem[Urquhart(2001)]{urquhart2001}
A.~Urquhart.
\newblock \emph{Basic Many-Valued Logic}, pages 249--295.
\newblock Springer Netherlands, Dordrecht, 2001.
\newblock ISBN 978-94-017-0452-6.
\newblock \doi{10.1007/978-94-017-0452-6_4}.

\bibitem[W{\'o}jcicki(1988)]{W88}
R.~W{\'o}jcicki.
\newblock \emph{Theory of logical calculi. {B}asic theory of consequence operations}, volume 199 of \emph{Synthese Library}.
\newblock Reidel, Dordrecht, 1988.
\newblock URL \url{http://books.google.com/books?id=Ucxfry7oowIC}.

\end{thebibliography}

\end{document}